\documentclass[12pt,a4paper]{article}
\usepackage{amscd,amsmath,amssymb,amsfonts,times}
\usepackage{mathrsfs}
\usepackage[all]{xy}
\numberwithin{equation}{section}
   \newtheorem{thm}{Theorem}[section]
   \newtheorem{lem}[thm]{Lemma}
   \newtheorem{sublem}[thm]{Sublemma}
   \newtheorem{fact}[thm]{Fact}
   \newtheorem{set}[thm]{Setting}
   
   \newtheorem{cor}[thm]{Corollary}

   \newtheorem{defn}[thm]{Definition}
   \newtheorem{exmp}[thm]{Example}
   \newtheorem{rem}[thm]{Remark}
\newcommand{\qed}
{\mbox{}\nolinebreak$\square$\medbreak\par}
\newenvironment{pf}{\par\smallskip\noindent\emph{Proof.}}{\hfill\qed\par\medskip}
\newenvironment{pf*}[1]{\par\smallskip\noindent\emph{#1.}}{\hfill\qed\par\medskip}
\begin{document}
\title{On $p$-adic vanishing cycles of log smooth families}
\author{Shuji Saito\footnote{Supported by Grant-in-Aid for Scientific research 15H03606, 17K18723}\; and Kanetomo Sato\footnote{Supported by Grant-in-Aid for Scientific research 16K05072}}
\date{}
\maketitle

%
%
%
%
%
\def\abs{{\sf abs}}
\def\cd{{\sf cd}}
\def\codim{{\sf codim}}
\def\ch{{\sf ch}}
\def\CH{{\sf CH}}
\def\cl{{\sf cl}}
\def\Coker{{\sf Coker}}
\def\cone{{\sf Cone}}
\def\dlog{{d{\text{\rm log}}}}
\def\dR{{\sf dR}}
\def\dX{d}
\def\dY{N}
\def\et{{\sf \acute{e}t}}
\def\sEnd{{\mathscr E}\hspace{-1.5pt}nd}
\def\sExt{{\mathscr E}\hspace{-1.5pt}xt}
\def\ffe{{e}}
\def\Frac{{\sf Frac}}
\def\gp{{\sf gp}}
\def\gr{\text{\rm gr}}
\def\gys{{\sf Gys}}
\def\h{{\sf h}}
\def\H{H}
\def\Hom{{\sf Hom}}
\def\sHom{{\mathscr H}\hspace{-2pt}om}
\def\id{{\sf id}}
\def\Image{{\sf Im}}
\def\ker{{\sf Ker}}
\def\lin{{\sf lin}}
\def\loc{{\sf loc}}
\def\log{\text{\rm log}}
\def\morph{\theta}
\def\Ob{{\sf Ob}}
\def\ord{{\sf ord}}
\def\Proj{{\sf Proj \,}}
\def\qis{\text{\rm qis}}
\def\red{\text{\rm red}}
\def\Res{{\sf Res}}
\def\rmF{\text{\rm F}}
\def\sgn{{\sf sgn}}
\def\sh{\,{\sf sh}}
\def\sing{{\sf sing}}
\def\Spec{\text{\rm Spec}}
\def\ssm{\smallsetminus}
\def\triv{\text{\rm triv}}
\def\U{\mathscr U}
\def\Uz{{\mathscr U}^0\hspace{-1.2pt}}
\def\Ua{{\mathscr U}^1\hspace{-1.6pt}}
\def\Ub#1{{\mathscr U}^{#1}\hspace{-1.4pt}}
\def\V{\mathscr V}
\def\val{\text{\rm val}}
\def\zar{\text{\rm Zar}}
\def\bA{\mathbb A}
\def\bC{\mathbb C}
\def\bF{\mathbb F}
\def\bK{\mathbb K}
\def\bL{\mathbb L}
\def\bN{\mathbb N}
\def\bP{\mathbb P}
\def\bQ{\mathbb Q}
\def\bR{\mathbb R}
\def\bT{\mathbb T}
\def\bZ{\mathbb Z}

\def\cB{\mathscr B}
\def\cC{\mathscr C}
\def\cD{\mathscr D}
\def\cE{\mathscr E}
\def\cF{\mathscr F}
\def\cH{\mathscr H}
\def\cI{\mathscr I}
\def\cK{\mathscr K}
\def\cL{\mathscr L}
\def\cM{\mathscr M}
\def\cN{\mathscr N}
\def\cO{\mathscr O}
\def\cU{\mathscr U}
\def\cV{\mathscr V}
\def\cY{\mathscr Y}
\def\cZ{\mathscr Z}
\def\ccS{\mathscr V}

\def\sfS{\text {\it Shv}}

\def\fD{\mathfrak D}
\def\fE{\mathfrak E}
\def\fF{\mathfrak F}
\def\fL{\mathfrak L}
\def\fM{\mathfrak M}
\def\fO{\mathfrak O}
\def\fT{\mathfrak T}
\def\fY{\mathfrak Y}
\def\fb{\mathfrak b}
\def\fd{\mathfrak d}
\def\fe{\mathfrak e}
\def\fk{\mathfrak k}
\def\fl{\mathfrak l}
\def\fm{\mathfrak m}
\def\fn{\mathfrak n}
\def\fp{\mathfrak p}

\def\tK{\text{\it K}}
\def\tM{\text{\it M}}

%
\def\ep{\epsilon}
\def\ka{\kappa}
\def\Lam{\varLambda}
\def\lam{\lambda}
\def\Pii{\xi}
\def\vG{\varGamma}
\def\ze{\zeta}
%
%
%
%
\def\ra{\rightarrow}
\def\lra{\longrightarrow}
\def\Lra{\Longrightarrow}
\def\la{\leftarrow}
\def\lla{\longleftarrow}
\def\Lla{\Longleftarrow}
\def\da{\downarrow}
\def\hra{\hookrightarrow}
\def\lmt{\longmapsto}
\def\sm{\setminus}
\def\bs#1{\boldsymbol{#1}}
\def\wt#1{\widetilde{#1}}
\def\wh#1{\widehat{#1}}
\def\spt{\sptilde}
\def\ol#1{\overline{#1}}
\def\ul#1{\underline{#1}}
\def\us#1#2{\underset{#1}{#2}}
\def\os#1#2{\overset{#1}{#2}}
\def\lim#1{\us{#1}{\varinjlim}}
\def\tcB{\wt{\cB}}
\def\tcZ{\wt{\cZ}}
\def\tom{\wt{\omega}}
\def\tomx#1#2{W_{\hspace{-2pt}#1}{\hspace{1pt}}\tom_{x,\log}^{#2}}
\def\tomY#1{\wt{\omega}^{#1}_Y}
\def\tomYZ#1{\wt{\omega}^{#1}_{Y^\circ}}
\def\tthY{\wt{\cD}_Y}
\def\tthYZ{\wt{\cD}_{Y^\circ}}
\def\thetalog{\theta^{\log}}
\def\tSigma{\wt{\varSigma}}
\def\qc{*\hspace{-0.6pt}}
\def\iim{*\hspace{-1.0pt}}
\def\fcI{f^*\hspace{-2.8pt}\cI}
\def\iI{i^*\hspace{-2.8pt}\cI}
\def\iL{i^*\hspace{-1.8pt}\cL}
\def\ilL{L} 

%
%
%
\def\Gm{{\mathbb G}_{\hspace{-1pt}\text{m}}}
\def\Ga{{\mathbb G}_{\hspace{-1pt}\text{a}}}
\def\zpn{{\bZ/p^n}}

\def\isom{\os{\simeq}{\lra}}
\def\lisom{\os{\simeq}{\lla}}
%
%
%
%
\def\mwitt#1#2#3{W_{\hspace{-2pt}#2}{\hspace{1pt}}\omega_{#1}^{#3}}
\def\witt#1#2#3{W_{\hspace{-2pt}#2}{\hspace{1pt}}\Omega_{#1}^{#3}}
\def\mlogwitt#1#2#3{W_{\hspace{-2pt}#2}{\hspace{1pt}}\omega_{{#1},{\log}}^{#3}}
\def\tlogwitt#1#2#3{W_{\hspace{-2pt}#2}{\hspace{1pt}}\tom_{{#1},{\log}}^{#3}}
\def\logwitt#1#2#3{W_{\hspace{-2pt}#2}{\hspace{1pt}}\Omega_{{#1},{\log}}^{#3}}
\def\loglogwitt#1#2{W_{\hspace{-2pt}n}{\hspace{1pt}}\Omega^{#2}_{(Y_{#1}\hspace{-1.5pt},\hspace{1pt}E_{#1})\hspace{-1pt},\hspace{1pt}\log}}
%
%
\section{Introduction}\label{sect1}

Let $K$ be a henselian discrete valuation field of mixed characteristic $(0,p)$, with residue field $k$.
Let $O_K$ be the ring of integers in $K$, and let $X$ be a regular scheme which is flat of finite type over $\Spec(O_K)$.
We consider cartesian squares of schemes
\[ \xymatrix{ X_K \; \ar[d] \ar@{}[rd]|{\square} \ar@<-1pt>@{^{(}->}[r]^-j & X \ar[d] \ar@{}[rd]|{\square} & \ar@<1pt>@{_{(}->}[l]_-i \ar[d]  \; X_k \\
 \Spec(K) \; \ar@<-1pt>@{^{(}->}[r] & \Spec(O_K) & \ar@<1pt>@{_{(}->}[l] \; \Spec(k).} \]
The Kummer short exact sequence of \'etale sheaves on $X_K$
\[ 0 \lra \mu_{p^n} \lra \cO_{X_K}^\times \os{\times p^n}\lra \cO_{X_K}^\times \lra 0 \]
yields a long exact sequence of \'etale sheaves on $X_k$
\[ 0 \lra i^*j_*\mu_{p^n} \lra i^*j_*\cO_{X_K}^\times \os{\times p^n}\lra i^*j_*\cO_{X_K}^\times \os{\delta}\lra 
i^*R^1j_*\mu_{p^n} \lra i^*R^1j_*\cO_{X_K}^\times \lra \dotsb. \]
Since $X$ is regular, we have $i^*R^1j_*\cO_{X_K}^\times =0$ and the connecting map $\delta$ in this sequence induces an isomorphism
\[ i^*R^1j_*\mu_{p^n} \cong \Coker\Big( i^*j_*\cO_{X_K}^\times \os{\times p^n}\lra i^*j_*\cO_{X_K}^\times\Big). \]
A motivation of this note is to extend this fundamental fact to higher cohomological degrees.
More precisely, we are concerned with the surjectivity of a geometric version of Tate's norm residue homomorphism
\[ \varrho^q_{X,n} : \cK^\tM_q/p^n \lra i^*R^qj_*\mu_{p^n}^{\otimes q}, \]
where $\cK^\tM_q$ denotes a Milnor \tK-sheaf defined as a quotient of $(i^*j_*\cO_{X_K}^\times)^{\otimes q}$
and $\cK^\tM_q/p^n$ denotes the cokernel of the multiplication by $p^n$ on $\cK^\tM_q$,
 cf.\ \S\ref{sect4} below.
The sheaf $i^*R^qj_*\mu_{p^n}^{\otimes q}$ on the right hand side is, so called, {\it the sheaf of $p$-adic vanishing cycles}, which is an \'etale sheaf of arithmetic and geometric interest.
In their paper \cite{BK}, Bloch and Kato proved that the map $\varrho^q_{X,n}$ is surjective in the case where $X$ is smooth over $\Spec(O_K)$. Later in his paper \cite{H}, Hyodo extended this surjectivity to the case where $X$ is a semistable family over $\Spec(O_K)$.
These surjectivity facts play a fundamental role in a construction of $p$-adic period maps in the $p$-adic Hodge theory cf.\ \cite{Ku}, \cite{K1}, \cite{K4}, \cite{T1}, \cite{T2}, \cite{YY}.

To state our main results more precisely, we introduce the following generalized situation with log poles.
Let $D$ be a normal crossing divisor on $X$ which is flat over $\Spec(O_K)$, and let 
\[ \xymatrix{ \psi : U:=X - (X_k \cup D) \; \ar@<-1pt>@{^{(}->}[r] & X }\]
be the natural open immersion. We then have a version of symbol map with log poles
\[ \varrho^q_{(X,D),n} : \cK^\tM_q/p^n \lra M_n^q:=i^*R^q\psi_*\mu_{p^n}^{\otimes q}, \]
where $\cK^\tM_q$ is again a Milnor \tK-sheaf defined as a quotient of $(i^*\psi_*\cO_U^\times)^{\otimes q}$,
 cf.\ \S\ref{sect4} below.
Now we state a main result of this paper, where quasi-log smoothness is a generalization of log-smoothness
 (cf.\ Definition \ref{cond5-1}, Example \ref{ex:logsmooth}):

\begin{thm}[Theorem \ref{thm5-1}]\label{thm1-1}
If $(X,D)$ is quasi-log smooth over $\Spec(O_K)$ and $K$ contains a primitive $p$-th root of unity, then $\varrho^q_{(X,D),n}$ is surjective for any $q \geq 2$.
\end{thm}

In his paper \cite{T2}, Tsuji proves an isomorphism in the derived category of \'etale sheaves on $X_k$
\stepcounter{equation}
\begin{equation}\label{eq1-1}
  {\mathcal S}_n(q)_{(X,D)} \cong \tau_{\leq q} i^*R\psi_*\mu_{p^n}^{\otimes q}
\end{equation}
for quasi-log smooth $(X,D)$ assuming $0 \leq q \leq p-2$,
where ${\mathcal S}_n(q)_{(X,D)}$ denotes a log syntomic complex.
His strategy is to show that both hand sides in \eqref{eq1-1} are invariant under log blow-ups, and then to reduce his assertion to the case that $X$ is smooth over $\Spec(O_K)$ (and $D = \emptyset$). This last case is due to Kurihara \cite{Ku}.
We will prove Theorem \ref{thm1-1} using his arguments on log blow-ups, which is the first key ingredient of our results.
We have to note that Theorem \ref{thm1-1} does not follow from \eqref{eq1-1}.
Indeed, it is not clear that the $q$-th cohomology sheaf of ${\mathcal S}_n(q)_{(X,D)}$ is generated by symbols,
 which is rather a consequence of \eqref{eq1-1} and Theorem \ref{thm1-1}.
\par
To continue the outline of our proof of Theorem \ref{thm1-1}, we introduce a subsheaf $\Ua \cK^\tM_q$ of $\cK^\tM_q$,
which is the subsheaf generated by the image of $i^*(1+\cI)^\times \otimes (i^*\psi_*\cO_U^\times)^{\otimes (q-1)}$,
where $\cI$ denotes the ideal sheaf of $\cO_X$ defining the reduced part $Y:=(X_k)_\red$ of $X_k$,
 and $(1+\cI)^\times$ means the kernel of the map $\cO_X^\times \to i_*\cO_Y^\times$.
We will further introduce a multi-index descending filtration on $\Ua \cK^\tM_q$, where the multi-indexes are assigned
 to irreducible components of $X_k$, cf.\ Definition \ref{defn4-1}.
To investigate the map $\varrho^q_{(X,D),n}$, we will need to control the behavior of the sheaf
\[ \Ua M_1^q := \Image (\Ua \cK^\tM_q \to M_1^q) \]
and a certain absolute logarithmic differential sheaf $\tom_{Y,\log}^q$ under log blow-ups (see Lemma \ref{lem5-3} below),
which corresponds to a key computation by Hyodo in the semistable family case, cf.\ \cite{H} Lemma (3.5).
Our second key ingredient is the computations on the multi-graded quotients of the induced filtration on $\Ua M_1^q$,
 which will be carried out by ideas of Kato, who introduced a new Cartier operator on the absolute differential modules with log poles, cf.\ \S\S\ref{sect3}--\ref{sect4}.
By this computation on multi-graded quotients, the behavior of $\Ua M_1^q$ and $\tom_{Y,\log}^q$ will be calculated by
 standard facts on the vanishing of the higher direct image of the structure sheaf under log blow-ups.
We would like to mention also that the idea of our computation on multi-graded quotients of $\Ua \cK^\tM_q$ has been used in a recent joint paper of the first author with R\"ulling \cite{RS}.
\par
\bigskip
Throughout this paper, we will work with the setting \ref{set1} stated below.
\stepcounter{thm}
\begin{defn}\label{def1-1}
{\rm Let $k$ be a field.
\begin{enumerate}
\item[{\rm(1)}]
{\it A normal crossing varity} over $k$ is a pure-dimensional scheme $Y$ which is separated of finite type over $k$ and everywhere \'etale locally
 isomorphic to \[ \Spec\big(k[T_1,\dotsc,T_\dY]/(T_1\dotsb T_a)\big) \quad \hbox{for some $1 \le a \le \dY = \dim(Y)+1$}.\]
\item[{\rm(2)}]
{\it An admissible divisor} on a normal crossing varity $Y$ is a reduced effective Cartier divisor $E$
such that the immersion $E \hra Y$ is everywhere \'etale locally isomorphic to \[\hspace{-10pt}\xymatrix{ \Spec\big(k[T_1,\dotsc,T_\dY]/(T_1\dotsb T_a,T_{a+1} \dotsb T_{a+b})\big) \;\ar@<-1pt>@{^{(}->}[r] & \Spec\big(k[T_1,\dotsc,T_\dY]/(T_1\dotsb T_a)\big) }\] for some $a,b \ge 1$ with $a+b \le \dY$.
\end{enumerate}
}\end{defn}
\noindent
\par
Let $K$ be a henselian discrete valuation field of characteristic $0$ whose residue field $k$ has characteristic $p>0$. Let $O_K$ be the integer ring of $K$. Unless mentioned otherwise, we do {\it not} assume that $k$ is perfect. Put $B:=\Spec(O_K)$ and $s:=\Spec(k)$.
\begin{set}\label{set1} $X$ is a regular scheme of finite type over $B$, and $D$ is a reduced divisor on $X$ which is flat over $B$ $(D$ may be empty$)$. We put $Y:=(X \times_B s)_{\red}$ and $U:=X \ssm (Y \cup D)$, and assume the following two conditions{\rm:}
\begin{itemize}
\item {\it The divisor $Y \cup D$ has normal crossings on $X$.} \par
\item {\it $Y$ is a normal crossing variety over $s$, and $(D \times_B s)_\red$ is an admissible divisor on $Y$.}
\end{itemize}
When $k$ is perfect, the first condition implies the second condition.
\end{set}

\smallskip
\section{Absolute differential modules with log poles}\label{sect2}
\medskip
Let the notation be as in Setting \ref{set1}, and let $i$ and $\psi$ be as follows:
\[\xymatrix{ Y \; \ar@<-1pt>@{^{(}->}[r]^-i & X & \ar@<1pt>@{_{(}->}[l]_-{\psi} \; U=X\ssm(Y\cup D). }\]
Put \[ \cL:=\psi_*\cO_U^{\times} \cap \cO_X \subset \psi_*\cO_U, \] which we regard as a sheaf of commutative monoids by the multiplicaition of functions. Let \[ \alpha : \iL \lra \cO_Y \] be the natural map of \'etale sheaves, where $i^*$ denotes the topological inverse image of \'etale sheaves. In this section, we study the following \'etale sheaves.
\newpage
\begin{defn}\label{def2-2}
{\rm\begin{enumerate}
\item[{\rm(1)}]
Let $\Omega^1_{Y/\bZ}$ be the usual absolute K\"ahler differential sheaf on $Y_{\et}$. We define the \'etale sheaf $\tom_Y^1$ on $Y$ as the quotient sheaf of \[ \Omega^1_{Y/\bZ} \oplus \big(\cO_Y \otimes_{\bZ} i^\iim \psi_*\cO_U^{\times}\big)\] divided by the $\cO_Y$-submodule generated by local sections of the form \[ (d\alpha(x),0) - (0,\alpha(x) \otimes x) \;\;\hbox{with} \;\; x \in i^\iim \cL. \] There is a logarithmic differential map \[ \dlog: i^\iim \psi_*\cO_U^{\times} \lra \tom_Y^1, \quad x \longmapsto (0,1 \otimes x).\] Put $\tom_Y^0:=\cO_Y$ and $\tom_Y^q := \wedge^q_{\cO_Y}\,\tom_Y^1$ for $q \ge 2$.
\item[{\rm(2)}]
We define $\tcZ_Y^q$ {\rm (}resp.\ $\tcB_Y^q${\rm )} as the kernel of $d:\tom_Y^q \to \tom_Y^{q+1}$ {\rm (}resp.\ the image of $d:\tom_Y^{q-1} \to \tom_Y^q${\rm )}, and put
\[ \tom_{Y,\log}^q:= \Image \big(\dlog :(i^\iim \psi_*\cO_U^{\times})^{\otimes q} \to \tom_Y^q \big). \]
\end{enumerate}}
\end{defn}
\begin{rem}\label{rem2-1}
{\rm The natural map $\cL \to \cO_X$ gives a log structure on $X$ in the sense of {\rm\cite{K2}}. In terms of log schemes, the sheaf $\tom_Y^1$ means the differential module $\omega^1_{(Y,L)/\bZ}$ defined in loc.\ cit., {\rm (1.7)}, where $L$ denotes the inverse image log structure of $\cL$ onto $Y$ {\rm (}loc.\ cit., {\rm (1.4))}.}
\end{rem}
\begin{thm}\label{thm2-1}
\begin{enumerate}
\item[{\rm(1)}]
The sheaf $\tom_Y^q$ is locally free over $\cO_Y$.
\item[{\rm(2)}]
There is a unique isomorphism
\[ C^{-1}  : \tom_Y^q \isom \cH^q(\tom_Y^\bullet) = \tcZ_Y^q/\tcB_Y^q \]
sending a local section $x \cdot \dlog(y_1) \wedge \dotsb \wedge \dlog(y_q)$ with $x \in \cO_Y$ and each $y_i \in i^\iim \psi_*\cO_U^{\times}$, to $x^p \cdot \dlog(y_1) \wedge \dotsb \wedge \dlog(y_q)+\tcB_Y^q$.
\item[{\rm(3)}]
There is a short exact sequence on $Y_\et$
\[ \xymatrix{ 0 \ar[r] & \tom_{Y,\log}^q \ar[r] & \tcZ_Y^q \ar[rr]^{1-C^{-1}\quad} && \cH^q(\tom_Y^{\bullet}) \ar[r] & 0. } \]
\end{enumerate}
\end{thm}
\begin{pf}
We first reduce the problem to the case that $X$ is a regular semistable family over $B=\Spec(O_K)$. Note that we may work \'etale locally. Indeed, once we prove (2) \'etale locally, then the isomorphisms $C^{-1}$ patch together automatically by the uniqueness. Assume the following three conditions (cf.\ Setting \ref{set1}):
\begin{itemize}
\item {\it $X$ is affine, and $Y = \Spec \big(k[T_1,\dotsc,T_\dX]/(T_1\dotsb T_a) \big)$ for $1 \le {}^{\exists} a \le \dX=\dim(X)$.}
\item {\it The irreducible components of $D$ are regular and principal.}
\item {\it There exists a regular sequence $t_1,\dotsc,t_{\dX}$ of prime elements of $\vG(X,\cO_X)$ such that $t_{\lam}$ lifts $T_{\lam} \in \vG(Y,\cO_Y)$ for $1 \le {}^\forall \lam \le \dX$ and such that $t_{a+1},\dotsc, t_{a+b}$ are uniformizers of the irreducible components of $D$ for $0 \le {}^\exists b \le \dX-a$.}
\end{itemize}
Let $\pi$ be a prime element of $O_K$. We have
\addtocounter{equation}{3}
\begin{equation}\label{eq2-a} \pi=u \, t_1^{e_1}\,t_2^{e_2}\dotsb t_a^{e_a}  \end{equation} for some $u \in \vG(X,\cO_X^\times)$ and some $e_1,\dotsc,e_a \ge 1$.
Put \begin{align*} X' & := \Spec\big(O_K[S_1,\dotsb,S_\dX]/(S_1\dotsb S_a - \pi)\big), \\ Y' & :=(X')_s = \Spec\big(k[S_1,\dotsb,S_\dX]/(S_1\dotsb S_a)\big), \\ U' & :=  \Spec\big(K[S_1,\dotsb,S_\dX,S_{a+1}^{-1},\dotsc,S_{a+b}^{-1}]/(S_1\dotsb S_a - \pi)\big). \end{align*}
Let $\psi'$ (resp.\ $i'$) be the open (resp.\ closed) immersion $U' \hra X'$ (resp.\ $Y' \hra X'$), and let $\beta$ be the isomorphism of schemes
\[ \beta : Y \isom Y', \qquad S_\lam \mapsto T_\lam \;\; (1 \le \lam \le \dX). \]
Put $\cK := \ker(\cO_X^\times \to \cO_Y^\times)$ and $\cK' := \ker(\cO_{X'}^\times \to \cO_{Y'}^\times)$.
 By \eqref{eq2-a}, $\psi_*\cO_U^\times/\cO_X^{\times}$ is a free abelian sheaf generated by $t_1,\dotsc,t_{a+b}$. Similarly $\psi'_*\cO_{U'}^\times/\cO_{X'}^{\times}$ is a free abelian sheaf generated by $S_1,\dotsc,S_{a+b}$. Hence there is an isomorphism of sheaves on $Y_{\et}$
\[ \beta^\iim \big(i'{}^\iim \psi'_*\cO_{U'}^\times/\cK'\big) \isom i^\iim \psi_*\cO_U^\times/\cK \]
that extends the isomorphism $\beta^\qc \cO_{Y'}^{\times} \isom \cO_Y^\times$ and sends $S_\lam \mapsto t_\lam$ for $1 \le \lam \le a+b$. By this isomorphism we see that $\beta^* \tom_{Y'}^q \cong \tom_Y^q$. Thus we are reduced to the case that $X$ is a regular semistable family over $B$.\par
We assume that $X$ is a regular semistable family over $B$ in what follows. Let $\omega_Y^q$ be the cokernel of the map
\[ \tom_Y^{q-1} \lra \tom_Y^q, \;\; x \mapsto \dlog(\pi) \wedge x. \]
We have a short exact sequence of complexes
\begin{equation}\label{eq2-0} \xymatrix{ 0  \ar[r] & \omega_Y^{\bullet}[-1] \ar[rr]^{\dlog(\pi) \wedge} && \tom_Y^\bullet \ar[r] &\omega_Y^\bullet \ar[r] & 0 } \end{equation}
and a short exact sequence of the $q$-th cohomology sheaves for any $q \ge 0$
\begin{equation}\label{eq2-1} \xymatrix{ 0 \ar[r] & \cH^{q-1}(\omega_Y^{\bullet}) \ar[rr]^{\dlog(\pi) \wedge} && \cH^q(\tom_Y^\bullet) \ar[r] &\cH^q(\omega_Y^\bullet) \ar[r] & 0,} \end{equation}
cf.\ \cite{T2} Lemma A.7 with $m=0$. We recall here the following facts due to Tsuji, loc.\ cit.\ Theorems A.3 and A.4 (cf.\ \cite{K2} Proposition (3.10), Theorem (4.12)\,(1)):
\addtocounter{thm}{3}
\begin{fact}\label{fact2-1}
\begin{enumerate}
\item[{\rm(1)}]
$\omega_Y^q$ is a locally free $\cO_Y$-module, and there is a unique isomorphism
\begin{equation}\notag
  C^{-1} : \omega_Y^q \isom \cH^q(\omega_Y^\bullet)
\end{equation}
sending a local section of the form $x \cdot \dlog(y_1) \wedge \dotsb \wedge \dlog(y_q)$ with $x \in \cO_Y$ and each $y_i \in i^\iim \psi_*\cO_U^{\times}$, to $x^p \cdot \dlog(y_1) \wedge \dotsb \wedge \dlog(y_q)+d\omega_Y^{q-1}$.
\item[{\rm(2)}]
Let $V$ be an open subset of $Y$ which is smooth over $k$ and for which $D \times_X V$ is empty. Then the short exact sequence \eqref{eq2-0} splits on $V$ as complexes. Consequently the exact sequence \eqref{eq2-1} splits on $V$, i.e., we have \[ \cH^q(\tom_V^\bullet) \isom \cH^q(\Omega_V^{\bullet}) \oplus \cH^{q-1}(\Omega_V^{\bullet}).\]
\item[{\rm(3)}]
There is a short exact sequence on $Y_\et$ \[\xymatrix{ 0 \ar[r] & \omega_{Y,\log}^q \ar[r] & \cZ_Y^q \ar[rr]^{1-C^{-1}\quad} && \cH^q(\omega_Y^\bullet) \ar[r] & 0, }\] where $\omega_{Y,\log}^q$ is defined as $\Image \big( \dlog :(i^\iim \psi_*\cO_U^{\times})^{\otimes q} \to \omega_Y^q \big)$.
\end{enumerate}
\end{fact}
Theorem \ref{thm2-1}\,(1) follows from \eqref{eq2-0} and Fact \ref{fact2-1}\,(1). We prove Theorem \ref{thm2-1}\,(2). Let $V$ be a dense open subset of $Y$ which is smooth over $k$ and for which $D \times_X V$ is empty. Let $\sigma$ be the open immersion $V \hra Y$. We first show that the canonical adjunction map
\stepcounter{equation}
\begin{equation}\label{eq2-3}
\cH^q(\tom_Y^\bullet) \lra \sigma_*\sigma^\iim \cH^q(\tom_Y^\bullet) \cong \sigma_*\big(\cH^q(\Omega_V^{\bullet}) \oplus \cH^{q-1}(\Omega_V^{\bullet})\big)
\end{equation}
is injective. Indeed by \eqref{eq2-1} there is a commutative diagram with exact rows
\[\xymatrix{ 0 \ar[r] & \cH^{q-1}(\omega_Y^{\bullet}) \ar[rr]^{\dlog(\pi) \wedge} \ar[d]& & \cH^q(\tom_Y^\bullet) \ar[r] \ar[d]_{\eqref{eq2-3}} & \cH^q(\omega_Y^\bullet) \ar[r] \ar[d] & 0 \\ 0 \ar[r] & \sigma_*\cH^{q-1}(\Omega_V^{\bullet}) \ar[rr]^{\dlog(\pi) \wedge} & & \sigma_*\cH^q(\tom_V^\bullet) \ar[r] & \sigma_* \cH^q(\Omega_V^\bullet) \ar[r] & 0, }\]
where the vertical arrows are adjunction maps. The left and right vertical arrows are injective by Fact \ref{fact2-1}\,(1). Hence the map \eqref{eq2-3} is injective. We define the map $C^{-1}:\tom_Y^q \to \cH^q(\tom_Y^\bullet)$ as follows. Using differential symbols, we easily see that the image of the composite map
\[\xymatrix{ \tom_Y^q \ar[r] & \sigma_*\big(\Omega_V^q \oplus \Omega_V^{q-1}\big) \ar[rr]^{C^{-1}\qquad\;\;} && \sigma_*\big(\cH^q(\Omega_V^{\bullet}) \oplus \cH^{q-1}(\Omega_V^{\bullet})\big). }\]
is contained in $\cH^q(\tom_Y^\bullet)$. We thus obtain the map $C^{-1}:\tom_Y^q \to \cH^q(\tom_Y^\bullet)$. By the construction, $C^{-1}$ is described by the local assignment as in Theorem \ref{thm2-1}\,(2), which implies the uniqueness of $C^{-1}$. Moreover it is bijective by the following commutative diagram with exact rows and Fact \ref{fact2-1}\,(1):
\[\xymatrix{ 0 \ar[r]  & \omega_Y^{q-1} \ar[rr]^{\dlog(\pi) \wedge} \ar[d]_{C^{-1}}^{\hspace{-1pt}\wr} && \tom_Y^q \ar[r]  \ar[d]_{C^{-1}} & \omega_Y^q \ar[r]  \ar[d]_{C^{-1}}^{\hspace{-1pt}\wr} & 0 \\ 0 \ar[r] & \cH^{q-1}(\omega_Y^{\bullet}) \ar[rr]^{\dlog(\pi) \wedge} & & \cH^q(\tom_Y^\bullet) \ar[r] & \cH^q(\omega_Y^\bullet) \ar[r] & 0. }\] This completes the proof of Theorem \ref{thm2-1}\,(2). \par
We prove Theorem \ref{thm2-1}\,(3). By the local presentation of $C^{-1}$, it is easy to see that the map $1-C^{-1}:\tcZ_Y^q \to \cH^q(\tom_Y^\bullet)$ is surjective and that its kernel contains $\tom_{Y,\log}^q$. Put
\[ L:=\ker\big(1-C^{-1}:\tcZ_Y^q \to \cH^q(\tom_Y^\bullet)\big). \]
We show that the natural inclusion $\tom_{Y,\log}^q \hra L$ is surjective. By \eqref{eq2-0}, \eqref{eq2-1} and Fact \ref{fact2-1}\,(1), there is a commutative diagram with exact rows
\[\xymatrix{ 0 \ar[r] & \cZ_Y^{q-1} \ar[rr]^{\dlog(\pi) \wedge} \ar[d]_{1-C^{-1}} && \tcZ^q_Y \ar[r] \ar[d]_{1-C^{-1}} & \cZ^q_Y \ar[r] \ar[d]_{1-C^{-1}} & 0  \\ 0 \ar[r] & \cH^{q-1}(\omega_Y^{\bullet}) \ar[rr]^{\dlog(\pi) \wedge} && \cH^q(\tom_Y^{\bullet}) \ar[r] & \cH^q(\omega_Y^{\bullet}) \ar[r] & 0, }\] By this diagram and Fact \ref{fact2-1}\,(3), the lower row of the following commutative diagram of complexes is exact:
\[ \xymatrix{ 0 \ar[r] & \omega_{Y,\log}^{q-1} \ar[rr]^{\dlog(\pi) \wedge} \ar@{=}[d] && \tom_{Y,\log}^q \ar[r] \ar[d] & \omega_{Y,\log}^q \ar[r] \ar@{=}[d] & 0 \\ 0 \ar[r] & \omega_{Y,\log}^{q-1} \ar[rr]^{\dlog(\pi) \wedge} && L \ar[r] & \omega_{Y,\log}^q \ar[r] & 0. } \]
Hence the middle vertical arrow is surjective, and we obtain Theorem \ref{thm2-1}\,(3).
\end{pf}

\smallskip
\section{Another Cartier isomorphism}\label{sect3}
\medskip
Let the notation be as in Setting \ref{set1}. Let $i$ and $\psi$ be as follows: \[\xymatrix{ Y \; \ar@<-1pt>@{^{(}->}[r]^-i & X & \ar@<1pt>@{_{(}->}[l]_-\psi \; U=X\ssm(Y\cup D). }\]
Let $\{Y_\lam\}_{\lam \in \Lam}$ be the irreducible components of $Y$.
For $\lam \in \Lam$, let $\cI_\lam \subset \cO_X$ be the defining ideal of $Y_\lam$. For $\fm=(m_\lam)_{\lam \in \Lam} \in \bN^\Lam$, put \[\cI^{(\fm)} := \prod_{\lam \in \Lam}\; \cI_\lam^{m_\lam} \subset \cO_X, \]
where $\bN$ denotes the set of natural numbers $\{0,1,2,\dotsc \}$.
\begin{defn}
{\rm\begin{enumerate}
\item[{\rm(1)}]
We endow $\bN^\Lam$ with a semi-order as follows. For $\fm=(m_\lam)_{\lam \in \Lam}$ and $\fn=(n_\lam)_{\lam \in \Lam} \in \bN^\Lam$, we say that $\fm \le \fn$ if $m_\lam \le n_\lam$ for all $\lam \in \Lam$.
\item[{\rm(2)}] We put ${\bf 0} :=(0)_{\lam \in \Lam}$ and ${\bf 1} :=(1)_{\lam \in \Lam} \in \bN^\Lam$.
\item[{\rm(3)}]
For $\fm,\fl \in \bN^\Lam$, we define a sheaf $\omega_{\fm,\fl}^q$ $(q \ge 0)$ on $Y_{\et}$ as
\[ \omega_{\fm,\fl}^q := \cI^{(\fm)}/\cI^{(\fm+\fl)} \otimes_{\cO_X} \tom_Y^q, \]
which is a locally free module over $\cO_X/\cI^{(\fl)}$ if $\fl \leqq {\bf 1}$.
We define a map $d_\fm^q:\cI^{(\fm)} \otimes_{\cO_X}\tom_Y^q \to \cI^{(\fm)} \otimes_{\cO_X}\tom_Y^{q+1}$ by the local assignment
\[\prod_{\lam \in \Lam} \, \pi_\lam^{m_\lam} \otimes \omega \; \mapsto \;
 \prod_{\lam \in \Lam}\, \pi_\lam^{m_\lam} \otimes \left(d\omega + \sum_{\lam \in \Lam} m_\lam\cdot\dlog(\pi_\lam) \wedge \omega \right) \quad (\omega \in \tom_Y^q),\]
where $\pi_\lam \in \cO_X$ denotes a local uniformizer of $Y_\lam$ for each $\lam \in \Lam$. This map does not depend on the choice of local uniformizers $\{\pi_\lam \}_{\lam \in \Lam}$. One can easily check that $d_\fm^{q+1} \circ d_\fm^q=0$ and that $d_\fn^q$ is compatible with $d_{\fm}^q$ for $\fn \ge \fm$. Hence $d_\fm^q$ induces a differential operator \[ d=d_{\fm,\fl}^q : \omega_{\fm,\fl}^q \lra \omega_{\fm,\fl}^{q+1}. \] Using this $d$, we regard $\omega_{\fm,\fl}^\bullet=(\omega_{\fm,\fl}^\bullet,d)$ as a complex.
\item[{\rm(4)}]
We define $\cZ_{\fm,\fl}^q$ $($resp.\ $\cB_{\fm,\fl}^q)$ as the kernel of $d:\omega_{\fm,\fl}^q \to \omega_{\fm,\fl}^{q+1}$ $($resp.\ the image of $d:\omega_{\fm,\fl}^{q-1} \to \omega_{\fm,\fl}^q)$, which are \'etale subsheaves of $\omega_{\fm,\fl}^q$.
\end{enumerate}}
\end{defn}
The following result is due to Kato:
\begin{thm}\label{thm3-1}
Let $\fm=(m_\lam)_{\lam \in \Lam}$ and $\fl=(\ell_\lam)_{\lam \in \Lam}$ be elements of $\bN^\Lam$ with ${\bf 0} \le \fl \le {\bf 1}$. Let $\fm' \in \bN^\Lam$ be the smallest element that satisfies $p\cdot \fm' \ge \fm$, and define $\fl'=(\ell'_\lam)_{\lam \in \Lam}$ by \[ \ell'_\lam := \begin{cases} 1 \quad & \hbox{$($if $\ell_\lam=1$ and $p|m_\lam)$} \\ 0 & (\hbox{otherwise}). \end{cases} \] Then there is an isomorphism \[ C^{-1}:\omega_{\fm',\fl'}^q \isom \cH^q(\omega_{\fm,\fl}^{\bullet}) = \cZ_{\fm,\fl}^q/\cB_{\fm,\fl}^q, \quad x \otimes \omega \mapsto x^p \otimes \omega, \] where $x$ $($resp.\ $\omega)$ denotes a local section of $\cI^{(\fm')}$ $($resp.\ $\tom_{Y,\log}^q)$.
\end{thm}
We first state an immediate consequence of this theorem. For $\mu \in \Lam$,
put $\tom_{Y_\mu}^q := \cO_{Y_\mu}\otimes_{\cO_Y}\tom_Y^q$ and define $\fd_{\mu}=(\delta_{\mu,\lam})_{\lam \in \Lam}$ by \[\delta_{\mu,\lam} := \begin{cases} 1 \quad & (\lam=\mu) \\ 0 & (\lam \ne \mu). \end{cases}\] Then we have $\omega_{\fm,\fd_\mu}^q = \cI^{(\fm)} \otimes_{\cO_X} \tom_{Y_\mu}^q$, and Theorem {\rm \ref{thm3-1}} implies the following{\rm :}
\begin{cor}\label{rem3-2}
\begin{enumerate}
\item[{\rm(1)}] The complex $\omega_{\fm,\fd_\mu}^{\bullet}$ is acyclic $($i.e., exact$)$ if $p \hspace{-3pt} \not | \, m_\mu$.
\item[{\rm(2)}] If $p|m_\mu$, then we have an isomorphism \[C^{-1}:\omega_{\fm',\fd_\mu}^q \isom \cH^q(\omega_{\fm,\fd_\mu}^{\bullet}),\] where $\fm' \in \bN^\Lam$ is the smallest element satisfying $p\cdot \fm' \ge \fm$.
\end{enumerate}
\end{cor}
\begin{pf*}{\it Proof of Theorem \ref{thm3-1}}
Let $\pi_{\lam} \in \cO_X$ be a local uniformizer of $Y_\lam$ for each $\lam \in \Lam$. If $p$ divides $m_\lam$ for any $\lam \in \Lam$, then the map $d:\omega_{\fm,\fl}^q \to \omega_{\fm,\fl}^{q+1}$ sends
\[\prod_{\lam \in \Lam} \, \pi_\lam^{m_\lam} \otimes \omega \;
  \mapsto \; \prod_{\lam \in \Lam} \, \pi_\lam^{m_\lam} \otimes d\omega,\]
and the assertion follows from Theorem \ref{thm2-1}\,(2). \par We prove the general case. Take a sequence of elements of $\bN^\Lam$ \[ \fm=\fm_0 \le \fm_1 \le \fm_2 \le \dotsb \le \fm_t=p\cdot\fm'\] such that \[\sum_{\lam \in \Lam} \ m_{i+1,\lam} - \sum_{\lam \in \Lam} \ m_{i,\lam} = 1 \quad \hbox{ for \, $0 \le i < t$},\] where $\fm_i=(m_{i,\lam})_{\lam \in \Lam}$ and $\fm_{i+1}=(m_{i+1,\lam})_{\lam \in \Lam}$. For $0 \le i \le t$, define $\fl_i=(\ell_{i,\lam})_{\lam \in \Lam}$ by \[ \ell_{i,\lam} := \begin{cases} 1 \quad & \hbox{(if $\ell_\lam=1$ and $m_{i,\lam}=m_\lam)$} \\ 0 & (\hbox{otherwise}). \end{cases} \] We have inclusions of complexes \[\omega_{\fm,\fl}^\bullet= \omega_{\fm_0,\fl_0}^\bullet \supset \omega_{\fm_1,\fl_1}^\bullet \supset \omega_{\fm_2,\fl_2}^\bullet \supset \dotsb \supset \omega_{\fm_t,\fl_t}^\bullet = \omega_{p\cdot\fm',\fl'}^\bullet, \] and exact sequences $(0 \le i < t)$ \[ 0 \lra \omega_{\fm_{i+1},\fl_{i+1}}^\bullet \lra \omega_{\fm_i,\fl_i}^\bullet \lra \omega_{\fm_i,\fd_\mu}^\bullet/\cI^{(\fl_i)}\omega_{\fm_i,\fd_\mu}^\bullet \lra 0, \] where $\mu=\mu(i)$ is the unique element of $\Lam$ such that $m_{i+1,\mu} > m_{i,\mu}$ and $\fd_\mu$ is as we defined before Corollary \ref{rem3-2}. It is enough to show the following two assertions:
\begin{enumerate}
\item[(A)]
{\it The complex $\omega_{\fm_i,\fd_\mu}^{\bullet}/\cI^{(\fl_i)}\omega_{\fm_i,\fd_\mu}^{\bullet}$ is acyclic for $0 \le i <t$.}
\item[(B)]
{\it There is an isomorphism $C^{-1}: \omega_{\fm',\fl'}^q \isom \cH^q(\omega_{p \cdot \fm',\fl'}^{\bullet})$.} 
\end{enumerate}
The assertion (B) follows from the proved case of the theorem. Because $\cI^{(\fl_i)}$ is locally free over $\cO_X$ for any $0 \le i < t$, the assertion (A) is reduced to the following:
\begin{lem}\label{lem3-3}
Let $\mu \in \Lam$ and assume $p \hspace{-3pt} \not | \, m_\mu$. Then the complex $\omega_{\fm,\fd_\mu}^{\bullet}$ is acyclic.
\end{lem}
\noindent
We prove this lemma in what follows. Let $\tom_{Y_\mu}^q$ be as before Corollary \ref{rem3-2}. Note that $\tom_{Y_\mu}^q$ is generated by $\Omega_{Y_\mu}^q$ (usual K\"ahler $q$-forms) and $q$-forms of the form $\dlog(\pi_\lam) \wedge \eta$ with $\lam \in \Lam$ and $\eta \in \Omega_{Y_\mu}^{q-1}$.
For $q \ge 1$, there is a residue homomorphism \[\Res^q : \tom_{Y_\mu}^q \lra \tom_{Y_\mu}^{q-1} \] characterized by the following two properties:
\begin{enumerate}
\item[(1)]
{\it For $\omega \in \Omega_{Y_\mu}^q$, $\Res^q(\omega)$ is zero.}
\item[(2)]
{\it For $\eta \in \Omega_{Y_\mu}^{q-1}$, we have \[\Res^q(\dlog(\pi_\lam) \wedge \eta)= \begin{cases}\eta \quad & (\lam=\mu) \\ 0 & (\lam \ne \mu). \end{cases}\]}
\end{enumerate}
Since $\omega_{\fm,\fd_\mu}^q = \cI^{(\fm)} \otimes_{\cO_X} \tom_{Y_\mu}^q$, we define a residue homomorphism \[\Res^q : \omega_{\fm,\fd_\mu}^q \lra \omega_{\fm,\fd_\mu}^{q-1} \] by $\Res^q(a \otimes \omega):=a \otimes \Res^q(\omega)$ for $a \in \cI^{(\fm)}$ and $\omega \in \tom_{Y_\mu}^q$. We show that
\addtocounter{equation}{4}
\begin{equation}\label{eq3-1}
d\,\Res^q(x) + \Res^{q+1}(dx)=m_\mu \cdot x \quad \hbox{ for any }\, x \in \omega_{\fm,\fd_\mu}^q,
\end{equation}
which implies that $\omega_{\fm,\fd_\mu}^{\bullet}$ is acyclic if $p \hspace{-3pt} \not | \, m_\mu$. Put \[ \Pii^\fm:=\prod_{\lam \in \Lam} \pi_\lam^{m_\lam} \in \prod_{\lam \in \Lam} \cI_\lam^{m_\lam}= \cI^{(\fm)}. \] For $x=\Pii^\fm \otimes \omega$ with $\omega \in \Omega_{Y_\mu}^q$, we have $\Res^q(x)=0$ and \begin{align*}\Res^{q+1}(dx) & = {\textstyle \Res^{q+1} \big(\Pii^\fm \otimes \big(d\omega + \sum_{\lam \in \Lam} \ m_\lam\cdot\dlog(\pi_\lam) \wedge \omega \big)\big)} \\ & = \Pii^\fm \otimes m_\mu \cdot \omega = m_\mu \cdot x. \end{align*}
For $x=\Pii^\fm\otimes \dlog(\pi_\nu) \wedge \eta$ with $\eta \in \Omega_{Y_\mu}^{q-1}$ and $\nu\ne\mu$, we have $\Res^q(x)=0$ and \begin{align*} \Res^{q+1}(dx) & = {\textstyle \Res^{q+1} \big(\Pii^\fm \otimes \big(-\dlog(\pi_\nu) \wedge d\eta + \sum_{\lam \in \Lam} \ m_\lam\cdot\dlog(\pi_\lam) \wedge \dlog(\pi_\nu) \wedge \eta \big)\big)} \\ & = {\textstyle \Pii^\fm \otimes m_\mu \cdot \dlog(\pi_\nu) \wedge \eta = m_\mu \cdot x.} \end{align*}
Finally for $x=\Pii^\fm\otimes \dlog(\pi_\mu) \wedge \eta$ with $\eta \in \Omega_{Y_\mu}^{q-1}$, we have \[{\textstyle d\,\Res^q(x) = d\big(\Pii^\fm \otimes \eta\big) = \Pii^\fm \otimes \left(d\eta + \sum_{\lam \in \Lam} \ m_\lam\cdot\dlog(\pi_\lam) \wedge \eta \right) }\] and \begin{align*} \Res^{q+1}(dx) &= {\textstyle \Res^{q+1} \big(\Pii^\fm \otimes \big(-\dlog(\pi_\mu) \wedge d\eta + \sum_{\lam \in \Lam} \ m_\lam\cdot\dlog(\pi_\lam) \wedge \dlog(\pi_\mu) \wedge \eta \big)\big)} \\ &= {\textstyle \Pii^\fm \otimes \big(-d\eta - \sum_{\lam \ne \mu} \ m_\lam \cdot \dlog(\pi_\lam) \wedge \eta\big)=-d\,\Res^q(x)+m_\mu \cdot x.} \end{align*} Thus we obtain \eqref{eq3-1}, Lemma \ref{lem3-3} and Theorem \ref{thm3-1}.
\end{pf*}
\stepcounter{thm}
\begin{cor}\label{cor3-4}
$\cZ_{\fm,\fl}^q$ is generated by local sections of the following forms{\rm:}
\begin{enumerate}
\item[{\rm(1)}]
$\prod_{\lam \in \Lam} \pi_\lam^{m_\lam} \otimes \big( d \eta + \sum_{\lam \in \Lam} \ m_\lam\cdot \dlog(\pi_\lam) \wedge \eta \big)$ with $\eta \in \tom_Y^{q-1}$, where $\pi_{\lam} \in \cO_X$ is a local uniformizer of $Y_\lam$ for each $\lam \in \Lam$.
\item[{\rm(2)}]
$x^p \otimes \omega$ with $x \in \cI^{(\fm')}$ and $\omega \in \tom_{Y,\log}^q$, where $\fm' \in \bN^\Lam$ is the smallest element satisfying $p\cdot \fm' \ge \fm$.
\end{enumerate}
\end{cor}

\smallskip
\section{Structure of the sheaf $\bs{\Ua M_1^q}$}\label{sect4}
\medskip
Let the notation be as in Setting \ref{set1}. Let $i$ and $\psi$ be as follows:
\[\xymatrix{ Y \; \ar@<-1pt>@{^{(}->}[r]^-i & X & \ar@<1pt>@{_{(}->}[l]_-{\psi} \; U=X\ssm(Y\cup D). }\] 
For $q \ge 0$ and $n \ge 1$, we define \'etale sheaves $M_n^q$ and $\cK^\tM_q$ on $Y$ as \[ M^q_n := i^\iim R^q\psi_* \mu_{p^n}^{\otimes q} \quad \hbox{ and } \quad  \cK^\tM_q := \begin{cases} \bZ & (q=0) \\ i^\iim \psi_*\cO_U^\times & (q=1) \\ ( i^\iim \psi_*\cO_U^\times )^{\otimes q}/J_q & (q \ge 2). \end{cases} \] Here $J_q$ for $q \ge 2$ denotes the subsheaf of $( i^\iim \psi_*\cO_U^\times )^{\otimes q}$ generated by local sections $x_1 \otimes \dotsb \otimes x_q$ such that $x_r + x_{r'}=0$ or $1$ for some $1 \le r < r' \le q$. There is a homomorphism of \'etale sheaves (cf.\ \cite{BK} 1.2)
\begin{equation}\label{symb}
\varrho^q_{(X,D),n} : \cK^\tM_q/p^n \lra M_n^q,
\end{equation}
which is a geometric version of Tate's norm residue map. For $x_1,\dotsc,x_q \in i^\iim \psi_*\cO_U^\times$, we denote the image of $\{x_1,x_2,\dotsc,x_q \} \in \cK ^M_q$ under \eqref{symb} again by $\{x_1,x_2,\dotsc,x_q \}$.
\stepcounter{thm}
\begin{defn}\label{defn4-1}
{\rm\begin{enumerate}
\item[{\rm(1)}]
We define $\Uz \cK^\tM_q$ as the full-sheaf $\cK^\tM_q$, and $\Ua \cK^\tM_q$ as the subsheaf generated locally by symbols of the form \[ \{1+x,y_1,\dotsc,y_{q-1}\}  \;\; \hbox{ with $x \in \iI$ and $y_1,\dotsc,y_{q-1} \in i^\iim \psi_*\cO_U^\times$}, \] where $\cI \subset \cO_X$ denotes the defining ideal of $Y$. We define $\Uz M^q_n$ and $\Ua M^q_n$ as the image of $\Uz \cK^\tM_q$ and $\Ua \cK^\tM_q$ under the map \eqref{symb}, respectively.
\item[{\rm(2)}]
For $\fm=(m_\lam)_{\lam \in \Lam} \in \bN^\Lam$ with $\fm \ge 1$, we define $\Ub {(\fm)}\cK^\tM_q$ as the subsheaf generated locally by symbols of the form \[ \{1+x,y_1,\dotsc,y_{q-1}\}  \;\; \hbox{ with $x \in \iI^{(\fm)}$ and $y_1,\dotsc,y_{q-1} \in i^\iim \psi_*\cO_U^\times$}, \] where $\cI^{(\fm)}$ is as we defined in the previous section. We define $\Ub {(\fm)}M^q_n$ as the image of $\Ub {(\fm)}\cK^\tM_q$ under the map \eqref{symb}.
\item[{\rm(3)}]
We define $\fe=(e_\lam)_{\lam \in \Lam} \in \bN^\Lam$ as follows. For $\lam \in \Lam$, let $e_\lam$ be the absolute ramification index of the discrete valuation ring $\cO_{X,\eta_\lam}$, where $\eta_\lam$ denotes the generic point of $Y_\lam$. We put $e'_\lam:=pe_\lam/(p-1)$ for $\lam \in \Lam$ and $\fe':=(e'_\lam)_{\lam \in \Lam} \in \bQ^\Lam$.
\item[{\rm(4)}]
For $\fm=(m_\lam)_{\lam \in \Lam}$ and $\fn=(n_\lam)_{\lam \in \Lam} \in \bQ^\Lam$, we say that $\fm < \fn$ $($resp.\ $\fm \le \fn)$ if $m_\lam < n_\lam$ $($resp.\ $m_\lam \le n_\lam)$ for any $\lam \in \Lam$. \end{enumerate}}
\end{defn}
The following lemma is straight-forward, and left to the reader:
\begin{lem}\label{lem4-2}
\begin{enumerate}
\item[{\rm(1)}] We have $\Ua M^q_n=\Ub {({\bf 1})}M^q_n$, where ${\bf 1}$ denotes $(1)_{\lam \in \Lam} \in \bN^\Lam$.
\item[{\rm(2)}] $M_n^q$, $\Uz M_n^q$ and $\Ua M_n^q$ are contravariantly functorial in the pair $(X,D)$.
\end{enumerate}
\end{lem}
\par
The main result of this section is the following, which is also due to Kato:
\begin{thm}\label{thm4-1}
Let $\fm$ and $\fl$ be elements of $\bN^\Lam$.
\begin{enumerate}
\item[{\rm(1)}]
Assume ${\bf 1} \le \fm < \fe'+{\bf 1}$ and ${\bf 0} \le \fl \le {\bf 1}$. For each $\lam \in \Lam$, assume $\ell_\lam=0$ if $m_\lam \ge e'_\lam$. Then there is an isomorphism \[ \omega_{\fm,\fl}^{q-1}\big/\cZ_{\fm,\fl}^{q-1} \isom \Ub {(\fm)}M^q_1\big/\Ub {(\fm+\fl)}M^q_1 \] given by the local assignment \[ x \otimes \dlog(y_1) \wedge \dotsb \wedge \dlog(y_{q-1}) \mapsto \left\{ 1+\wt{x},y_1,\dotsc,y_{q-1} \right\}+\Ub {(\fm+\fl)}M^q_1 \] for $x \in i^\iim \big(\cI^{(\fm)}\big/\cI^{(\fm+\fl)}\big)$ and $y_1,\dotsc,y_{q-1}\in i^\iim \psi_*\cO_U^\times$, where $\wt{x}$ is a lift of $x$ to $\iI^{(\fm)}$.
\item[{\rm(2)}]
If $\fm \ge \fe'$, then $\Ub {(\fm)}M^q_1$ is zero.
\end{enumerate}
\end{thm}
Theorem {\rm \ref{thm4-1}} describes the structure of the sheaf $\Ua M^q_1$ as follows.
\begin{cor}\label{rem4-2}
Take a sequence of elements of $\bN^\Lam$ \[ {\bf 1} = \fm_0 \le \fm_1 \le \fm_2 \le \dotsb \le \fm_i \le \dotsb \le \fm_t \] satisfying the following conditions{\rm:}
\begin{enumerate}
\item[{\rm(a)}] $\fl_i:=\fm_{i+1}-\fm_i$ satisfies $\fl_i \le {\bf 1}$ for any $i \le t-1$.
\item[{\rm(b)}] We have $\fe' \le \fm_t < \fe' + 1$.\, {\rm(}Note that such $\fm_t$ is unique in $\bN^\Lam$.{\rm)}
\item[{\rm(c)}] For any $(i,\lambda)$ with $i \le t-1$ and $m_{i,\lam} \ge e'_\lam$, the $\lam$-component of $\fl_i$ is zero.
\end{enumerate}
We then have \begin{align*} \Ub{(\fm_i)}M^q_1\big/\Ub{(\fm_{i+1})}M^q_1 & \cong \omega_{\fm_i,\fl_i}^{q-1}\big/\cZ_{\fm_i,\fl_i}^{q-1} \quad \hbox{ for }\; 0 \le i \le t-1 \\ \Ub{(\fm_t)}M^q_1& =0 \end{align*} by Theorem {\rm \ref{thm4-1}\,(1)} and {\rm (2)}, respectively.
\end{cor}
\noindent
To prove Theorem \ref{thm4-1}, we need the following lemma:
\begin{lem}\label{lem4-3}
\begin{enumerate}
\item[{\rm(1)}]
For $\fm,\fn \in \bN^\Lam$, we have $\{\Ub {(\fm)}\cK^\tM_q,\Ub {(\fn)}\cK^\tM_{q'}\} \subset \Ub{(\fm+\fn)}\cK^\tM_{q+q'}$.
\item[{\rm(2)}]
There is a surjective homomorphism \[ i^\iim \cO_X \otimes (i^\iim \psi_*\cO_U^{\times})^{\otimes r} \lra \tom_Y^r, \quad x \otimes y_1 \otimes \dotsb \otimes y_r \mapsto \ol x \cdot \dlog(y_1) \wedge \dotsb \wedge \dlog(y_r),\] where for $x \in i^\iim \cO_X$, $\ol x$ denotes its residue class in $\cO_Y$. The kernel of this map is generated by local sections of the following forms{\rm:}\smallbreak
\begin{itemize} \item $x \otimes y_1 \otimes \dotsb \otimes y_r$ with $x \in \iI$ or $y_s \in i^\iim (1+\cI)$ for some $1 \le s \le r$, \item $x \otimes y_1 \otimes \dotsb \otimes y_r$ with $y_s=y_{s'}$ for some $1 \le s < s' \le r$, \item $\sum_{s=1}^{m}\, (x_s \otimes x_s \otimes y_1 \otimes \dotsb \otimes y_{r-1}) - \sum_{t=1}^{\ell}\, (x'_t \otimes x'_t \otimes y_1 \otimes \dotsb \otimes y_{m-1})$ such that all $x_s$ and $x'_t$ belong to $i^\iim (\cO_X \cap \psi_*\cO_U^{\times})$ and such that the sums ${\sum}_{s=1}^m \, x_s$ and ${\sum}_{t=1}^\ell \, x'_t$ taken in $i^\iim \cO_X$ satisfy ${\sum}_{s=1}^m\, x_s \equiv {\sum}_{t=1}^\ell\, x'_t \mod \iI$. \end{itemize} \end{enumerate}
\end{lem}
\begin{pf}
(1) follows from the same argument as in \cite{BK} Lemma 4.1. \par
(2) Let $z$ be a point on $Y$. Put $A:=\cO_{X,\ol z}^{\sh}$, $I:=\cI_{\ol z}$ and $L:=(\psi_*\cO_U^\times)_{\ol z}$. Let $A[L]$ be the free $A$-module over the set $L$. There is a surjective $A$-homomorphism $A[L] \to (\tom_Y^1)_{\ol z}$ sending $a[b] \mapsto \ol a \cdot \dlog(b)$. Its kernel is the $A$-submodule generated by elements of the following forms: \smallbreak
\begin{enumerate}
\item[(i)] $a[b]$ with $a \in I$ or $b \in 1+I$, \qquad \qquad (ii) $[b\cdot b']-[b]-[b']$ with $b,b' \in L$, \smallbreak
\item[(iii)] $\displaystyle \sum_{s=1}^m\, a_s[a_s]-\sum_{t=1}^\ell\, a'_t[a'_t]$ $(a_s,a'_t \in A \cap L)$ with $\displaystyle \sum_{s=1}^m \, a_s \equiv \sum_{t=1}^\ell \, a'_t \mod I$. \smallbreak
\end{enumerate}
The claim follows from this fact. The details are straight-forward and left to the reader.
\end{pf}
\begin{pf*}{\it Proof of Theorem \ref{thm4-1}}
(1) By Lemma \ref{lem4-3}, the local assignment in the theorem gives a well-defined surjective homomorphism of sheaves \[\xymatrix{ \rho_{\fm,\fl} : \omega_{\fm,\fl}^{q-1} \ar@{->>}[r] & \Ub{(\fm)}M^q_1\big/\Ub{(\fm+\fl)}M^q_1. }\] We prove $\rho_{\fm,\fl}\big(\cZ_{\fm,\fl}^{q-1}\big)=0$ assuming $\fm < \fe'$, locally on $Y$. We may assume that $Y_{\lam}$'s are principal on $X$. Fix uniformizers $\pi_{\lam} \in i^\iim \cO_X$ of $Y_\lam$ $(\lam \in \Lam)$ and put \[ \Pii^\fm:=\prod_{\lam \in \Lam} \pi_\lam^{m_\lam} \in \iI^{(\fm)}. \] It is enough to show that local sections of $\cZ_{\fm,\fl}^{q-1}$ of the forms (1) and (2) of Corollary \ref{cor3-4} map to zero under $\rho_{\fm,\fl}$.
For $y_1 \in i^\iim \cO_X^{\times}$ and $y_2,\dotsc,y_{q-1} \in i^\iim \psi_*\cO_U^{\times}$, we have \begin{align*} & \{1+\Pii^\fm y_1,y_1,y_2,\dotsc,y_{q-1} \} + \sum_{\lam \in \Lam} \ m_\lam \cdot \{1+\Pii^\fm y_1, \pi_\lam,y_2,\dotsc,y_{q-1}\} \\ & = \{1+\Pii^\fm y_1,\Pii^\fm y_1,y_2,\dotsc,y_{q-1}\}=-\{1+\Pii^\fm y_1,-1,y_2,\dotsc,y_{q-1}\} \in \Ub{\fm+\fl}\cK^\tM_q. \end{align*} Hence $\rho_{\fm,\fl}(\omega)=0$ for \begin{align*} \omega &=\Pii^\fm \otimes \big( d \eta + \textstyle \sum_{\lam \in \Lam} \ m_\lam\cdot \dlog(\pi_\lam) \wedge \eta \big) \in \cZ_{\fm,\fl}^{q-1} \\ &\hbox{ with } \;\; \eta = y_1 \cdot \dlog(y_2) \wedge \dotsb \wedge \dlog(y_{q-1}) \in \tom_Y^{q-2}. \end{align*} Next let $\fm' \in \bN^\Lam$ be the smallest element that satisfies $p\cdot \fm' \ge \fm$. For $x \in \iI^{(\fm')}$ and $y_1,\dotsc,y_{q-1} \in i^\iim \psi_*\cO_U^{\times}$, we have \[ \{1+x^p,y_1,\dotsc,y_{q-1} \} - p \cdot \{1+x,y_1,\dotsc,y_{q-1} \} \in \Ub{\fm'+\fe}\cK^\tM_q \subset \Ub{\fm+\fl}\cK^\tM_q, \] where we have used the assumption $\fm < \fe'$ to verify $\fm'+\fe \ge \fm+\fl$. Hence $\rho_{\fm,\fl}(\omega)=0$ for \[ \omega=x^p \otimes \dlog(y_1) \wedge \dotsb \wedge \dlog(y_{q-1}) \in \cZ_{\fm,\fl}^{q-1}.\] Thus we obtain  $\rho_{\fm,\fl}\big(\cZ_{\fm,\fl}^{q-1}\big)=0$ for $\fm < \fe'$.
\par We prove $\rho_{\fm,\fl}\big(\cZ_{\fm,\fl}^{q-1}\big)=0$ for $\fm < \fe'+{\bf 1}$. Define $\fn=(n_\lam)_{\lam \in \Lam}$ and $\fl'=(\ell'_\lam)_{\lam \in \Lam}$ as \[ n_\lam := \begin{cases}m_\lam \quad & (\hbox{if $m_\lam < e'_\lam$}) \\ m_\lam-1 \quad & (\hbox{if $m_\lam \ge e'_\lam$}) \end{cases} \quad \hbox{ and } \quad \ell'_\lam := \begin{cases}\ell_\lam \quad & (\hbox{if $m_\lam < e'_\lam$}) \\ 1 \quad & (\hbox{if $m_\lam \ge e'_\lam$}). \end{cases} \] We have $\fn \le \fm$, ${\bf 1} \le \fn < \fe'$, ${\bf 0} \le \fl' \le {\bf 1}$ and $\fn+\fl'=\fm+\fl$ by the assumptions on $\fm$ and $\fl$, and there is a commutative diagram \[\xymatrix{ \omega_{\fm,\fl}^{q-1} \; \ar@{^{(}->}[r] \ar[d]_{\rho_{\fm,\fl}} & \omega_{\fn,\fl'}^{q-1} \ar[d]^{\rho_{\fn,\fl'}} \\ \Ub{(\fm)}M^q_1/\Ub{(\fm+\fl)}M^q_1 \; \ar@{^{(}->}[r] & \Ub{(\fn)}M^q_1/\Ub{(\fn+\fl')}M^q_1,}\] where the top horizontal arrow maps $\cZ_{\fm,\fl}^{q-1}$ into $\cZ_{\fn,\fl'}^{q-1}$. Hence $\rho_{\fm,\fl}\big(\cZ_{\fm,\fl}^{q-1}\big)$ is zero by the previous case and the injectivity of the bottom horizontal arrow.

\par It remains to prove the injectivity of the induced map \[\ol \rho_{\fm,\fl} : \omega_{\fm,\fl}^{q-1}/\cZ_{\fm,\fl}^{q-1} \lra \Ub{(\fm)}M^q_1/\Ub{(\fm+\fl)}M^q_1. \] Since $\omega_{\fm,\fl}^{q-1}/\cZ_{\fm,\fl}^{q-1}$ is a subsheaf of $\omega_{\fm,\fl}^q$, the canonical adjunction map \[\begin{CD}\omega_{\fm,\fl}^{q-1}/\cZ_{\fm,\fl}^{q-1} \lra \bigoplus_{y \in Y^0}\ i_{y*}i_y^\iim \big(\omega_{\fm,\fl}^{q-1}\big/\cZ_{\fm,\fl}^{q-1}\big)\end{CD}\] is injective by Theorem \ref{thm2-1}\,(1), where for $y \in Y^0$, $i_y$ denotes the natural map $y \hra Y$. Hence we may replace $X$ with $\Spec(\cO_{X,y_\mu})$ ($\mu \in \Lam$), where $y_\mu$ denotes the generic point of $Y_\mu$. By the definition of $d:\omega_{\fm,\fl}^{q-1}\to \omega_{\fm,\fl}^q$, we have \[ \omega_{\fm,\fl}^{q-1}\big/\cZ_{\fm,\fl}^{q-1} \cong \begin{cases} \Omega_{y_\mu}^{q-1} \quad &(\hbox{if $p \hspace{-3pt} \not | \, m_\mu$, $\ell_\mu=1$ and $m_\mu < e'_\mu$}) \\ d\Omega_{y_\mu}^{q-1}\oplus d\Omega_{y_\mu}^{q-2} \quad &(\hbox{if $p | m_\mu$, $\ell_\mu=1$ and $m_\mu < e'_\mu$}) \\ 0 & (\hbox{otherwise}) \end{cases} \] and the assertion follows from \cite{BK} Corollary 1.4.1\,(ii)--(iv).
\par (2)
For $\fm \in \bN^\Lam$ with $\fm \ge \fe'$, $1+\cI^{(\fm)}$ is contained in $(1+\cI^{(\fm-\fe)})^p$. The assertion follows from this fact. \end{pf*}


\smallskip
\section{Surjectivity of the symbol map}\label{sect5}
\medskip
Let $O_K$ be as in \S\ref{sect1}, and let $\pi$ be a prime element of $O_K$. 
\begin{defn}\label{def5-1}
{\rm For an injective morphism of monoids
\[ h : \bN \lra \bN^\dX, \;\; 1 \mapsto (e_\lam)_{1 \le \lam \le d}, \]
we define a scheme $X^h$ and a divisor $D^h$ on $X^h$ as
\begin{align*}
X^h &
 := \Spec\Big(O_K\big[T_1,\dotsc,T_d \big]\big/
\big({\textstyle \prod_{\lam\,\text{ with }\,e_\lam \ge 1}} \ T_\lam^{e_\lam}-\pi\big)\Big) \\
 D^h & := \big\{{\textstyle\prod_{\lam\,\text{ with }\,e_\lam = 0}}\ T_{\lam}=0 \big\} \subset X^h.
\end{align*}
Put $Y^h:=(X^h)_{s,\red}$. We define a scheme $\cY^h$ as $\Spec(k[T_1,T_2,\dotsc,T_d])$ and denote the natural closed immersion $Y^h \hra \cY^h$ by $\iota^h$.
}\end{defn}
Let $(X,D)$ be as in Setting {\rm \ref{set1}}. We introduce the following terminology:
\begin{defn}\label{cond5-1}
{\rm
We say that $(X,D)$ is {\it quasi-log smooth over} $B=\Spec(O_K)$, if it is,
everywhere \'etale locally on $X$, isomorphic to $(X^h,D^h)$ for some injective morphism of monoids $h:\bN \to \bN^{\dX}$ and a prime element $\pi \in O_K$, where $\dX$ denotes $\dim(X)$.
}
\end{defn}
\begin{exmp}\label{ex:logsmooth}
{\rm
Let $(X,D)$ be a pair as in Setting \ref{set1}.
\begin{enumerate}
\item[(1)]
Let $\cM$ be the log structure on $X$ associated with $D$, and
let $\cN$ be the log structure on $B=\Spec(O_K)$ associated with the closed point $s \in B$.
If the canonical morphism $(X,\cM) \to (B,\cN)$ of log schemes is smooth in the sense of \cite{K2} (3.3), then the pair $(X,D)$
 is quasi-log smooth over $\Spec(O_K)$ in our sense. Note also that the converse is not necessarily true.
\item[(2)]
As a consequence of (1),
a pair $(X,D)$ as in Setting {\rm \ref{set1}} is quasi-log smooth over $\Spec(O_K)$, if the multiplicities of the irreducible components of $X_s$ are prime to $p$.
\end{enumerate}
}
\end{exmp}
\begin{thm}\label{thm5-1}
Assume that $K$ contains a primitive $p$-th root of unity $\zeta_p$,
 and that $(X,D)$ is quasi-log smooth over $B$. Then the symbol map in \eqref{symb}
\[ \varrho^q_{(X,D),1} : \cK^\tM_q/p^n \lra M_n^q \]
 is surjective, and there is an isomorphism
\[ M_1^q/\Ua M_1^q \isom \tom_{Y,\log}^q \] fitting into a commutative diagram
\addtocounter{equation}{4}
\begin{equation}\label{eq5-1}
\xymatrix{ & \cK^\tM_q/p \ar@{->>}[ld]_{\ol{\varrho^q_{(X,D),n}}} \ar@{->>}[rd]^{\dlog} & \\  M_1^q/\Ua M_1^q  \ar[rr]^-\simeq
 & & \tom_{Y,\log}^q. }
\end{equation}
\end{thm}
\noindent
Our proof of Theorem \ref{thm5-1} will be complete in the next section.
In this section, we reduce the theorem to Lemma \ref{lem5-3} below.
Let $y$ be a generic point of $Y$, and let $i_y$ be the natural map $y \hra Y$. The strict henselian local ring $\cO_{X,{\ol y}}^{\sh}$ is a discrete valuation ring by the regularity of $X$. Hence there is an isomorphism
\begin{equation}\label{eq5-2}
 i_y^\iim \big(M_1^q/\Ua M_1^q\big) \isom \Omega_{y,\log}^q \oplus \Omega_{y,\log}^{q-1} = i_y^\iim \,\tom_{Y,\log}^q
\end{equation}
(\cite{BK} Lemma 5.3).
Since $\tom_Y^q$ is locally free over $\cO_Y$ by Theorem \ref{thm2-1}\,(1), the adjunction map
\[ \tom_{Y,\log}^q \lra \bigoplus_{y\in Y^0}\ i_{y*}i_y^\iim \,\tom_{Y,\log}^q\]
is injective, and the isomorphism \eqref{eq5-2} induces a surjective map
\[ \xymatrix{ \Uz M_1^q/\Ua M_1^q \ar@{->>}[r]  & \tom_{Y,\log}^q } \]
(cf.\ \cite{S3} Lemma 2.3). This map fits into the diagram \eqref{eq5-1} with $M_1^q/\Ua M_1^q$ replaced by $\Uz M_1^q/\Ua M_1^q$.
In what follows, put
\[
 M^q:=M_1^q, \quad
 N^q := \ker \big(\Uz M^q/\Ua M^q \to \tom_{Y,\log}^q \big) \quad \hbox{ and } \quad L^q := M^q/\Uz M^q. \]
We have show that $N^q=0$ and $L^q=0$.
Note also that once we show $L^q=0$, we will have shown that $M_n^q=\Uz M_n^q$ for all $n \ge 1$ by a standard argument as in \cite{BK} Corollary 6.1.1.
We first prove that $N^q=L^q=0$ in a simple case:\par

\addtocounter{thm}{2}
\begin{lem}\label{lem5-1}
Assume that $(X,D)$ is quasi-log smooth over $B$ and that the underlying scheme $X$ is smooth over $B$.
Then we have $N^q=0$ and $L^q=0$, i.e., Theorem {\rm \ref{thm5-1}} holds for $(X,D)$.
\end{lem}
\begin{pf}
When $D=\emptyset$, the assertion follows from a theorem of Bloch-Kato \cite{BK} Theorem 1.4.
We proceed the proof of the lemma by induction on the number of the irreducible components of $D$.
Since the problem is \'etale local on $X$ and $X$ is smooth over $B$ by assumption,
we may suppose that $X=\Spec(O_K[T_2,T_3,\dotsc,T_d])$ and that $D=\{T_2 T_3 \dotsb T_r=0\} \subset X$,
where $2 \leqq r \leqq d=\dim(X)$.
\par
Fix an irreducible component $V$ of $D$, which is also smooth over $B$.
Put $D':=D - V$ as an effective Cartier divisor, and let $E$ be the pullback of $D'$ onto $V$.
It is easy to see that $E$ is a simple normal crossing divisor on $V$ and the pair $(V,E)$ is quasi-log smooth over $B$ and that the underlying scheme $V$ is smooth over $B$.
Recall that $Y:=X_{s,\red}\,(=X_s)$ and $U:=X\ssm(Y\cup D)$.
Now put $Z:=V_{s,\red}\,(=V_s)$, $U':=X\ssm(Y \cup D')$ and $W:=V\ssm(Z \cup E)$, and consider a commutative diagram of schemes
\[\xymatrix{
  Z \; \ar@<-1pt>@{^{(}->}[r]^{i'} \ar[d]_{\iota} & V \ar[d] & \ar@<1pt>@{_{(}->}[l]_-{\theta} \; W \ar[d]^\alpha \\
  Y \; \ar@<-1pt>@{^{(}->}[r]^i  & X  & \ar@<1pt>@{_{(}->}[l]_-{\psi'} \; U'
  & \ar@<1pt>@{_{(}->}[l]_\beta \ar@/^7mm/[ll]_\psi \; U. }\]
We then have a commutative diagram of \'etale sheaves on $Y$ whose upper row is a complex and whose lower row is exact
\addtocounter{equation}{1}
\begin{equation}\label{eq5-1-1}
\xymatrix{
\cK^\tM_{q,(X,D')}/p \ar[r] \ar[d]_{\varrho^q_1} & \cK^\tM_q/p \ar@{->>}[r] \ar[d]_{\varrho^q_1} & \iota_*\cK^\tM_{q-1,(V,E)}/p \ar[d]_{\varrho^{q-1}_1} \\
M_{(X,D')}^q \ar[r] & M^q \ar[r] & \iota_*M_{(V,E)}^{q-1}.
}\end{equation}
Here we put $M_{(X,D')}^q:=i^*R^q\psi'_*\mu_p^{\otimes q}$ and
$M_{(V,E)}^{q-1}:=i'^*R^{q-1}\theta_*\mu_p^{\otimes (q-1)}$, and the sheaves
$\cK^\tM_{q,(X,D')}$ and $\cK^\tM_{q-1,(V,E)}$ are Milnor \tK-sheaves defined as quotients of
 $(i^*\psi'_*\cO_{U'}^\times)^{\otimes q}$ and $(i'^*\theta_*\cO_W^\times)^{\otimes (q-1)}$, respectively.
The right arrow in the upper row is a boundary map of Milnor \tK-sheaves, which one can check to be surjective.
The lower row is obtained by applying $i^*R^q\psi'_*$ to the Gysin distinguished triangle on $(U')_\et$
\[ \alpha_*\mu_p^{\otimes (q-1)}[-2] \lra \mu_{p}^{\otimes q} \lra R\beta_*\mu_p^{\otimes q}
 \lra \alpha_*\mu_p^{\otimes (q-1)}[-1]. \]
Now the left and the right vertical arrows in \eqref{eq5-1-1}
 are surjective by the induction hypothesis, and we obtain the surjectivity of $\varrho^q_{(X,D),1}$ by a simple diagram chase,
 that is, $L^q=0$.
\par
To prove that $N^q=0$, we consider the following commutative diagram of sheaves on $Y_\et$, whose middle row is exact and whose other rows and columns are complexes:
\[\xymatrix{
  & \Ua M_{(X,D')}^q \ar[r] \ar[d] & \Ua M^q \ar[r]^-{r_1} \ar[d] & \iota_*\Ua M_{(V,E)}^{q-1} \ar[r] \ar[d] & 0 \\
 &  M_{(X,D')}^q \ar[r] \ar[d] & M^q \ar[r]^-{r_2} \ar[d] & \iota_*M_{(V,E)}^{q-1} \ar[r] \ar[d] & 0 \\
0 \ar[r] &  \tom_{Y',\log}^q  \ar[r] &  \tom_{Y,\log}^q  \ar[r]^-{r_3} & \iota_* \tom_{Z,\log}^{q-1}  \ar[r] & 0,
}\]
where $\tom_{Y',\log}^q$ (resp.\ $\tom_{Z,\log}^{q-1}$) denotes the logarithmic differential sheaf defined for the pair $(X,D')$ (resp.\ $(V,E)$); the arrow $r_1$ denotes the map induced by $r_2$, and $r_3$ denotes a residue map of logarithmic differential sheaves. 
By the smoothness of $X$ over $B$, we have
\[ \tom_{Y,\log}^q \cong \Omega^q_Y(\log \hspace{1pt} D_Y)_{\log} \oplus \Omega^{q-1}_Y(\log \hspace{1pt} D_Y)_{\log} \qquad \hbox{($D_Y := D \cap Y$)} \]
and similar presentations for $\tom_{Y',\log}^q$ and $\tom_{Z,\log}^{q-1}$.
By this fact and a simple variant of the purity of logarithmic differential sheaves \cite{Sh} Theorem 3.2, we see that the bottom row is exact.
On the other hand, one can easily check that the sequence
\[ \xymatrix{(0 \ar[r] &)\, \cU^m M_{(X,D')}^q \ar[r] & \cU^m M^q \ar[r]^-{r_1} & \iota_*\cU^m M_{(V,E)}^{q-1} \ar[r] & 0} \]
is exact for $1 \le m \le e':=pe/(p-1)$ by Theorem \ref{thm4-1} and descending induction on $m$.
Hence the top row in the above diagram is also exact. Now the left and the right columns are exact by the induction hypothesis, and the middle column is exact as well by a simple diagram chase, which shows that $N^q=0$.
\end{pf}
\par
\medskip
In the rest of this section, we reduce the general case of Theorem {\rm \ref{thm5-1}} to Lemma \ref{lem5-3} below.
Fix an arbitrary point $x \in Y$.
We show the stalks $(N^q)_{\ol x}$ and $(L^q)_{\ol x}$ are zero by induction on $c:=\codim_Y(x)$. If $c=0$, then $(N^q)_{\ol x}=(L^q)_{\ol x}=0$ by \eqref{eq5-2}. In what follows, assume $c \ge 1$ and the following induction hypothesis:
\begin{enumerate}
\item[(1)]
{\it $(N^q)_{\ol y}=(L^q)_{\ol y}=0$ for any $q \ge 0$ and any $y \in Y$ of codimension $ \le c-1$.}
\end{enumerate}
 Since the problem is \'etale local, we may assume $(X,D)=(X^h,D^h)$ for an injective morphism $h : \bN \hra \bN^\dX$ of monoids ($\dX=\dim(X)$). Sorting the components of $\bN^\dX$ if necessarily, we assume the following two conditions:
\begin{enumerate}
\item[(2)]
{\it The first component of $h(1) \in \bN^\dX$ is non-zero.}
\item[(3)]
{\it The composite map \[ \xymatrix{ Y \; \ar@{^{(}->}[r]^{\iota^h \qquad \qquad \quad} & \cY = \Spec(k[T_1,\dotsc,T_d]) \ar[r] & \Spec(k[T_{c+2},\dotsc,T_d])}\] sends $x$ to the generic point of $\Spec(k[T_{c+2},\dotsc,T_d])$ {\rm(}cf.\ \cite{T2} Lemma {\rm5.3)}.}
\end{enumerate}
Following the idea of Tsuji in \cite{T2} Proof of Theorem 5.1, we decompose $h : \bN \to \bN^\dX$ into a sequence of morphisms of monoids
\[ \begin{CD} h \; : \; \bN @>{h^0}>> \bN^\dX @>{\ka^1}>> \bN^\dX @>{\ka^2}>> \dotsb @>{\ka^{r}}>> \bN^\dX \end{CD}\]
which satisfies the following two conditions:
\begin{enumerate}
\item[(4)]
{\it $h^0(1)=(e,\overbrace{0,\dotsc,0}^{c\text{-copies}},*,\dotsc,*)$ for some $e \ne 0$}
 ({\it cf.}\ (2)).
\item[(5)]
{\it For $1 \le t \le \dX$, let $\ep_t \in \bN^\dX$ be the element whose $t$-th component is $1$ and whose other components are $0$. Then for $1 \le \nu \le r$, $\ka^{\nu}$ sends $\ep_t$ $(1 \le t \le \dX)$ to \[ \begin{cases} \ep_t & (t \ne m) \\ \ep_m + \ep_n & (t = m) \end{cases} \; \hbox{ for some } \, m \ne n \; \hbox{ with }\; 1 \le m \le c+1,\; 1 \le n \le c+1. \]}
\end{enumerate}
\def\X#1{{}^{#1}\hspace{-1.7pt}X}
\def\Y#1{{}^{#1}Y}
\def\UU#1{{}^{#1}U}
\def\D#1{{}^{#1}\hspace{-1.6pt}D}
Put $h^{\nu}:=\kappa^{\nu}\kappa^{\nu-1}\dotsb\kappa^1h^0$ and $\X \nu :=X^{h^{\nu}}$, and let $f^{\nu}$ be the morphism induced by $\ka^\nu$:
\[\begin{CD} f^{\nu} : \X \nu \lra \X {\nu-1}. \end{CD}\]
We further fix some notation. Put $\Y {\nu}:=Y^{h^{\nu}}=(\X {\nu})_{s,\red}$, and let $x^{\nu} \in \Y {\nu}$ be the image of $x \in Y$ under the composite
\[ \begin{CD} Y=\Y r @>{g^r}>> \Y {r-1} @>{g^{r-1}}>> \dotsb @>{g^{\nu+1}}>> \Y \nu, \end{CD} \]
where $g^\nu : \Y \nu \to \Y {\nu-1}$ denotes the morphism induced by $f^{\nu}$. For $0 \le \nu \le r$, let $\sigma^\nu$ be the composite map
\[ \Y \nu \os{\iota^\nu}{\hra} \cY^{h^{\nu}}=\Spec(k[T_1,\dotsc,T_d]) \lra \Spec(k[T_{c+2},\dotsc,T_d]), \]
where $\iota^\nu$ denotes $\iota^{h_\nu}$.
Since $\sigma^\nu=\sigma^{\nu-1}g^\nu$ by (5), the point $\sigma^\nu(x^\nu)$ is the generic point of $\Spec(k[T_{c+2},\dotsc,T_d])$ for any $0 \le \nu \le r$ by (3). This implies the following:
\begin{enumerate}
\item[(6)]
{\it For any $0 \le \nu \le r$, $x^\nu$ has codimension $c$ on $\Y \nu$. Consequently, $x^{\nu}$ is a closed point of $(g^{\nu})^{-1}(x^{\nu-1})$.
}
\end{enumerate}
We also need the following fact (cf.\ \cite{T2} Lemmas 3.2 and 3.4):
\begin{enumerate}
\item[(7)]
{\it For $1 \le \nu \le r$, $f^{\nu}$ factors as
\[ \xymatrix{ \X \nu \; \ar@<-1pt>@{^{(}->}[r] & \ol{\X \nu} \ar[r]^{{\ol {f^{\nu}}}\;\;} & \X {\nu-1}, }\]
where the left arrow is an open immersion and $\ol {f^{\nu}}$ is the blow-up at the closed subscheme $\{T_m=T_n=0\} \subset \X {\nu-1}$.
The fibers of ${\ol {f^{\nu}}}$ have dimension at most one.}
\end{enumerate}
\par
Put \[ \D \nu := D^{h_\nu} \quad \hbox{ and } \quad \UU \nu := \X\nu \ssm (\Y \nu \cup \D \nu), \] and define the sheaves ${}_{\nu}M^q$, ${}_{\nu}L^q$ and ${}_{\nu}N^q$ on $\Y \nu_{\et}$ for the diagram \[\xymatrix{ \Y \nu \; \ar@<-1pt>@{^{(}->}[r] & \X\nu & \ar@<1pt>@{_{(}->}[l] \; \UU \nu }\] in the same way as for $M^q$, $L^q$ and $N^q$ on $Y_{\et}$, respectively. In what follows, we prove \[ ({}_{\nu}L^q)_{\ol{x^\nu}} = ({}_{\nu}N^q)_{\ol{x^\nu}} = 0 \] by induction on $0 \le \nu \le r$. We first note:
\addtocounter{thm}{1}
\begin{lem} We have $({}_0L^q)_{\ol{x^0}}=0$ and $({}_0N^q)_{\ol{x^0}}=0$. \end{lem}
\begin{pf}
By the assumption (4),
the pair $(\X 0,\D 0)$ is, \'etale locally around $x^0$, isomorphic to a quasi-log smooth pair $(X',D')$ over a henselian discrete valuation ring $A'$ of mixed characteristic such that $X'$ is smooth over $A'$, cf.\ \cite{T2} the first isomorphism in p.\ 559.
Hence the assertion follows from Lemma \ref{lem5-1}.
\end{pf}
Assume $\nu \ge 1$ and the following induction hypothesis:
\begin{enumerate}
\item[(8)]
{\it $({}_{\nu-1}L^q)_{\ol{x^{\nu-1}}}=0$ and $({}_{\nu-1}N^q)_{\ol{x^{\nu-1}}}=0$.}
\end{enumerate}
We change the notation slightly and put
\[ \begin{cases} X := \Spec\big(\cO^{\sh}_{\X {\nu-1},\,\ol{x^{\nu-1}}}\big) \\ Y := X_{s,\red}\\ D := \Spec\big(\cO^{\sh}_{\D {\nu-1},\,\ol{x^{\nu-1}}}\big) \\ U := X\ssm(Y \cup D) \end{cases} \quad \hbox{ and } \quad \begin{cases} X' := \ol {\X \nu} \times_{\X {\nu-1}} \Spec\big(\cO^{\sh}_{\X {\nu-1},\,\ol{x^{\nu-1}}}\big) \\ Y' := (X')_{s,\red} \\ D' := \ol {\D {\nu}} \times_{\ol {\X \nu}} X' \\ U' := X'\ssm(Y' \cup D') \end{cases} \]
for simplicity. Here $\ol {\D \nu}$ denotes the closure of $(\ol {f^\nu})^{-1}(\D \nu)_\red \ssm (\ol {\X{\nu}})_{s,\red} \subset \ol {\X \nu}$. Note that $\D \nu=\ol {\D \nu} \times_{\ol {\X \nu}} \X \nu$. Let $i' : Y' \hra X'$ and $\psi' : U' \hra X'$ be the canonical closed and open immersions, respectively. We define \'etale sheaves $M^q$, $L^q$ and $N^q$ on $Y'$ as
\begin{equation*}
M^q :=i'{}^\iim R^q\psi'_*\mu_p^{\otimes q}, \quad
N^q := \ker \big(\Uz M^q/\Ua M^q \twoheadrightarrow \tom_{Y',\log}^q \big),
 \quad L^q := M^q/\Uz M^q.
\end{equation*}
In view of (6), once we prove $N^q$ and $L^q$ are zero, we will finish the induction on $\nu$ and $c$. 
We will prove the following lemma in the next section:
\begin{lem}[{\bf cf.\ \cite{H} Lemma (3.5)}]\label{lem5-3}
For any $t \ge 0$, we have \begin{align*} \vG(Y,\tom_{Y,\log}^t) & \cong \vG(Y',\tom_{Y',\log}^t), \\ \H^1(Y',\Ua M^t) & =\H^1(Y',\tom_{Y',\log}^t)=0. \end{align*}
\end{lem}
We prove here that $N^q$ and $L^q$ are zero admitting this lemma.
Noting that $\mu_p \cong \bZ/p\bZ$ on $U'$ by the assumption on $K$, we compute the Leray spectral sequence \[ E_2^{a,b}=\H^a(Y',M^b) \Lra \H^{a+b}(U',\mu_p^{\otimes q}) \cong \H^{a+b}(U,\mu_p^{\otimes q}), \]
where we have used the proper base-change theorem (\cite{SGA4} XII.5.2) for the identification $\H^a(Y',M^b) \cong \H^a(X',R^b\psi'_*\mu_p^{\otimes q})$ and also used the fact that $\ol{f^{\nu}}$ induces an isomorphism $U' \cong U$. Since $\cd_p(Y') \le 1$ by (7), this spectral sequence yields a short exact sequence \[ 0 \lra \H^1(Y',M^{q-1}) \lra \H^q(U,\mu_p^{\otimes q}) \lra \vG(Y',M^q) \lra 0. \] Because $L^t$ and $N^t$ are skyscraper sheaves on $Y'$ for any $t \ge 0$ by the induction hypothesis (1) for $X'$, both $\H^1(Y',M^{q-1})$ and $\H^1(Y',\Uz M^q)$ are zero by Lemma \ref{lem5-3}. Hence there is a commutative diagram whose lower row is exact
\[ \xymatrix{ & \Uz \H^q(U,\mu_p^{\otimes q}) \ar@{=}[r] \ar[d] & \H^q(U,\mu_p^{\otimes q}) \ar[d]_{\wr\hspace{-1pt}} \\
0 \ar[r] & \vG(Y',\Uz M^q) \ar[r] & \vG(Y',M^q) \ar[r] & \vG(Y',L^q) \ar[r] & 0, } \]
where $\Ub {\bullet} \H^q(U,\mu_p^{\otimes q})$ ($\bullet=0,1$) denotes the filtration on the stalk of the sheaf of $p$-adic vanishing cycles on $Y$ (cf.\ Definition \ref{defn4-1}\,(1), Lemma \ref{lem4-2}\,(2)), and the upper equality follows from the induction hypothesis (8).
This diagram shows that the skyscraper sheaf $L^q$ is zero. We next show that $N^q$ is zero. Put \[ \gr_{\U}^0M^q:=\Uz M^q/\Ua M^q=M^q/\Ua M^q. \] Since $N^q$ is skyscraper by (1), there is an exact sequence
\[  0 \lra \vG(Y',N^q) \lra \vG(Y',\gr_{\U}^0M^q) \os{\alpha}{\lra} \vG(Y',\tom_{Y',\log}^q) \lra 0 \]
and a commutative diagram with exact rows (cf.\ Lemma \ref{lem4-2}\,(2))
\[ \xymatrix{ 0 \ar[r] & \Ua \H^q(U,\mu_p^{\otimes q}) \ar[r] \ar[d] & \H^q(U,\mu_p^{\otimes q}) \ar[r] \ar[d]_{\wr\hspace{-1pt}} & \vG(Y,\tom_{Y,\log}^q) \ar[r] \ar[d]_\beta & 0 \\ 0 \ar[r] & \vG(Y',\Ua M^q) \ar[r] & \vG(Y',M^q) \ar[r] & \vG(Y',\gr_{\U}^0M^q) \ar[r] & 0. }\]
Here the upper row is exact by the induction hypothesis (8), and the lower row is exact by Lemma \ref{lem5-3}. The arrow $\beta$ denotes the map induced by the left square, and the middle vertical arrow is bijective by the proof of the vanishing of $L^q$. Now this diagram shows that $\alpha$ is bijective, because $\beta$ is surjective and $\alpha\beta$ is bijective by Lemma \ref{lem5-3}. Hence the skyscraper sheaf $N^q$ is zero. Thus the induction on $\nu$ and $c$ is complete and we obtain Theorem \ref{thm5-1}, assuming Lemma \ref{lem5-3}.

\smallskip
\section{Proof of Lemma \ref{lem5-3}}\label{sect6}
\medskip
In this section we prove Lemma \ref{lem5-3} to finish the proof of Theorem \ref{thm5-1}.
Let the notation be as in Setting \ref{set1}. Assume that
\begin{align*}
 X & = \Spec\big(O_K [T_1,\dotsc,T_d]/(T_1^{e_1}\dotsb T_a^{e_a}-\pi)\big)
       \quad (e_1,\dotsc,e_a \ge 1,\; d=\dim(X)) \\
 D & = \{T_{a+1}\dotsb T_d=0 \} \subset X,
\end{align*}
where $\pi$ is a prime element of $O_K$ and $D$ is empty if $a=d$. Let \[ f : X' \lra X \] be the blow-up at the regular closed subscheme $\{T_b=T_c=0\} \subset X$ with $1 \le b < c \le d$.
Put $Y':=(X')_{s,\red}$. We define a reduced normal crossing divisor $D'$ on $X'$ as \[ D':=\ol { f^{-1}(D)_\red \ssm Y'}  \subset X'\] and define the sheaves $\tom_{Y',\log}^q$ and $M_{1,X'}^q$ on $(Y')_{\et}$ in the same way as for $\tom_{Y,\log}^q$ and $M_1^q$ on $Y_{\et}$. Let \[ g: Y' \lra Y \] be the morphism induced by $f$. Lemma \ref{lem5-3} follows from the following lemma:
\begin{lem}\label{lem6-1}
Let $D(Y_\et)$ be the derived category of complexes of \'etale sheaves on $Y$.
\begin{enumerate}
\item[{\rm(1)}] We have $\tom_{Y,\log}^q \isom Rg_*\tom_{Y',\log}^q$ in $D(Y_\et)$ for any $q \ge 0$.
\item[{\rm(2)}] We have $R^1g_*(\Ua M_{1,X'}^q)=0$ for any $q \ge 0$.
\end{enumerate}
\end{lem}

To prove this lemma, we introduce some notation, which will be useful throughout the proof of Lemma \ref{lem6-1}.
For $\lam \in \Lam:=\{ 1,2,\dotsc, a\}$, let $Y_\lam$ be the closed subset $\{T_\lam=0\} \subset X$ endowed with the reduced subscheme structure, which is an irreducible component of $Y$. Let $\{Y'_\lam \}_{\lam \in \Lam'}$ be the irreducible components of $Y'$. We have
\[
\Lam' =
\begin{cases}
 \Lam \cup \{ \infty \} \quad & (\hbox{if $b \le a$}) \\
 \Lam \quad & (\hbox{if $a < b$}),
\end{cases}
\]
where $Y'_\lam$ for $\lam \in \Lam$ means the strict transform of $Y_\lam$ and $Y'_{\infty}$ is the exceptional fiber of $f$. For $\fm=(m_\lam)_{\lam \in \Lam} \in \bN^\Lam$, we define $f^*\fm=(n_\lam)_{\lam \in \Lam'} \in \bN^{\Lam'}$ as \[ n_\lam = \begin{cases} m_\lam \quad & (\hbox{if $\lam \in \Lam$}) \\ m_b \quad & (\hbox{if $b \le a <c$ and $\lam = \infty$}) \\ m_b+m_c \quad & (\hbox{if $c \le a$ and $\lam = \infty$}). \end{cases} \] 
For $\lam \in \Lam'$, let $\cI'_\lam \subset \cO_{X'}$ be the defining ideal of $Y'_\lam$. For $\fn \in \bN^{\Lam'}$, we define $\cI_{X'}^{(\fn)} \subset \cO_{X'}$ in the same way as for $\cI^{(\bullet)} \subset \cO_X$, cf.\ \S\ref{sect3}.
We have
\stepcounter{equation}
\begin{equation}\label{eq6-0}
 \fcI^{(\fm)}=\cI_{X'}^{(f^\qc\fm)}.
\end{equation}
We will often write $\cI \subset \cO_X$ and $\cI' \subset \cO_{X'}$ for the defining ideals of $Y$ and $Y'$, respectively.
\stepcounter{thm}
\begin{sublem}\label{sublem6-1}
\begin{enumerate}
\item[{\rm(1)}]
We have $g^*\tom_Y^q \isom \tom_{Y'}^q$ for any $q \ge 0$, where $g^*$ denotes the inverse image of coherent sheaves.
\item[{\rm(2)}]
Assume that $c \le a$, and let $\fn=(n_\lam)_{\lam \in \Lam'} \in \bN^{\Lam'}$ be an arbitrary element with $n_b+n_c=n_\infty+1$.
Then we have $Rg_*(\cI_{X'}^{(\fn)}/\cI_{X'}^{(\fn+\fd_\infty)})=0$ in $D(Y_{\et})$.
See the definitions before Corollary {\rm\ref{rem3-2}} for $\fd_\infty \in \bN^{\Lam'}$.
\item[{\rm(3)}]
We have $\cO_Y \isom Rg_*\cO_{Y'}$ in $D(Y_{\et})$.
\end{enumerate}
\end{sublem}
\begin{pf*}{\it Proof of Sublemma \ref{sublem6-1}}\,
(1)
Let $(X,\cL)$ and $(X',\cL')$ be the log schemes associated with the pairs $(X,D)$ and $(X',D')$, respectively (cf.\ Remark \ref{rem2-1}), and let $L$ (resp.\ $L'$) be the inverse image log structure onto $Y$ (resp.\ onto $Y'$).
The morphism $(Y',L') \to (Y,L)$ of log schemes induced by $f$ is \'etale in the sense of \cite{K2} (3.3) by loc.\ cit.\ Theorem (3.5). Hence the assertion follows from loc.\ cit., Proposition 3.12.
\par
(2) When $c \le a$, $Y'_{\infty}$ is isomorphic to the trivial $\bP^1$-bundle over $Y_b \cap Y_c$, and $Y'_{\infty} \cap Y'_b$ and $Y'_{\infty} \cap Y'_c$ are relative hyperplane sections of $Y'_\infty$ over $Y_b \cap Y_c$.
Since $(\cI'_{\!\infty})^{n_\infty} \big/(\cI'_{\!\infty})^{n_\infty+1} \cong\cO(n_\infty)$ on $Y'_\infty$, we have
\[ \cI_{X'}^{(\fn)}\big/\cI_{X'}^{(\fn+\fd_\infty)} \cong (\cI'_b)^{n_b} \otimes_{\cO_{X'}} (\cI'_c)^{n_c} \otimes_{\cO_{X'}} \big((\cI'_{\!\infty})^{n_\infty} \big/(\cI'_{\!\infty})^{n_\infty+1} \big) \cong \cO(-1) \]
on $Y'_\infty$, where we have used the fact that $\cI'_\lam$ is principal for $\lam \ne b,c,\infty$
and the assumption that $n_b+n_c=n_\infty+1$.
Hence we obtain the assertion by a standard fact on the cohomology of projective lines (\cite{Ha} Chapter III Theorem 5.1).
\par
(3)
Noting that $f : X' \to X$ is the blow-up along a closed subscheme of $X$ defined by a regular sequence
and that $\cI=T_0T_1\dotsb T_a \cO_X$ is a free $\cO_X$-module (of rank $1$), we have
\[ \cO_X \isom Rf_*\cO_{X'} \quad \hbox{ and } \quad \cI \isom Rf_*\fcI  \quad \hbox{ in } \quad D(X_\et) \]
by the theorem on formal functions (\cite{Ha} Chapter III Theorem 11.1) and the cohomology of projective lines.
Our task is to check
\addtocounter{equation}{1}
\begin{equation}\label{eq6+1}
 Rf_*(\cI'/\fcI)=0 \quad \hbox{ in } \quad D(X_\et).
\end{equation}
Note that $\cI'=\cI_{X'}^{({\bf 1})}$ and that $\fcI=\cI_{X'}^{(f^{\qc}{\bf 1})}$ by \eqref{eq6-0}.
If $a<c$, then we have $f^*{\bf 1} = {\bf 1}$ in $\bN^{\Lam'}$ and the assertion \eqref{eq6+1} is obvious.
If $c \le a$, then we have
\begin{equation}\label{eq6++0}
  f^*{\bf 1}={\bf 1} + \fd_{\infty} \quad  \hbox{ in \;\; $\bN^{\Lam'}$}
\end{equation}
and the assertion follows from Sublemma \ref{sublem6-1}\,(2).
This completes the proof.
\end{pf*}

\begin{pf*}{\it Proof of Lemma \ref{lem6-1}}
(1) By Theorem \ref{thm2-1}\,(2) and (3), it is enough to show that
\[ {\mathrm {(a)}} \;\; \tom_Y^q \isom Rg_*\tom_{Y'}^q \quad \hbox{ and } \quad {\mathrm {(b)}} \;\; \tcZ_Y^q \isom Rg_*\tcZ_{Y'}^q  \quad \hbox{ in } \quad D(Y_\et).\]
We first show (a). Since $\tom_Y^q$ is locally free over $\cO_Y$ by Theorem \ref{thm2-1}\,(1), we have \[Rg_*\tom_{Y'}^q \cong Rg_*g^*\tom_Y^q \cong Rg_*(\cO_{Y'} \otimes_{\cO_{Y'}} g^*\tom_Y^q) \cong (Rg_*\cO_{Y'}) \otimes_{\cO_Y} \tom_Y^q \cong \tom_Y^q \]
by projection formula and Sublemma \ref{sublem6-1}\,(3).
We next show (b). The case $q=0$ follows from (a) and the isomorphism $\tcZ_Y^0=(\cO_Y)^p \cong \cO_Y$. We proceed the proof by induction on $q$. There is a commutative diagram with distinguished rows in $D(Y_\et)$ \[\xymatrix{ \tcZ_Y^q \ar[r] \ar[d] & \tom_Y^q \ar[r]^{d \;\;} \ar[d]_{\wr \hspace{-1pt}} & \tcB_Y^{q+1} \ar[r] \ar[d] & \tcZ_Y^q[1] \ar[d] \\ Rg_*\tcZ_{Y'}^q \ar[r] & Rg_*\tom_{Y'}^q \ar[r]^{d \;\;} & Rg_*\tcB_{Y'}^{q+1} \ar[r] & Rg_*\tcZ_{Y'}^q[1]. }\]
By Theorem \ref{thm2-1}\,(2), there is another commutative diagram with distinguished rows in $D(Y_\et)$ \[\xymatrix{ \tcB_Y^{q+1} \ar[r] \ar[d] & \tcZ_Y^{q+1} \ar[r]^{C} \ar[d] & \tom_Y^{q+1} \ar[r] \ar[d]_{\wr \hspace{-1pt}} & \tcB_Y^{q+1}[1] \ar[d] \\ Rg_*\tcB_{Y'}^{q+1} \ar[r] & Rg_*\tcZ_{Y'}^{q+1}\ar[r]^{C} & Rg_*\tom_{Y'}^{q+1} \ar[r] & Rg_*\tcB_{Y'}^{q+1}[1], } \] where $C$ denotes the inverse of the isomorphism $C^{-1}$. The induction on $q$ works by these diagrams.
\par
(2)
We first show
\begin{equation}\label{eq6-1}
R^1g_*\big(\Ub{(f^\qc{\bf 1})}M_{1,X'}^q\big)=0.
\end{equation}
Take a sequence of elements of $\bN^\Lam$ \[ {\bf 1} = \fm_0 \le \fm_1 \le \fm_2 \le \dotsb \le \fm_i \le \dotsb \le \fm_t \] satisfying the following conditions:
\begin{enumerate}
\item[(i)] If $i \le t-1$, then $\fl_i:=\fm_{i+1}-\fm_i$ agrees with $\fd_\mu$ for some $\mu=\mu(i) \in \Lam$. \item[(ii)] We have $\fm_i \ge \fe'$ (in $\bQ^\Lam$) if and only if $i=t$.
\end{enumerate}
See the definitions before Corollary \ref{rem3-2} for $\fd_\mu$, and see Definition \ref{defn4-1}\,(3) for $\fe' \in \bQ^\Lam$.
We define $\fE' \in \bQ^{\Lam'}$ for $X'$ in the same way as we defined $\fe'$ for $X$.
For $\fn,\fn' \in \bN^{\Lam'}$ with ${\bf 0} \le \fn' \le {\bf 1}$, we define the sheaf $\omega_{\fn,\fn',X'}^{q-1}$ on $Y'$ in the same way as for $\omega_{*,*}^{q-1}$ on $Y$. By the choice of the above $\fm_i$'s and the fact that $f^*\fe'=\fE'$, we see that
\begin{enumerate}
\item[(i$'$)] For $i \le t-1$, we have $f^*\fl_i = f^*\fm_{i+1}-f^*\fm_i \in \bN^{\Lam'}$ and satisfies ${\bf 0} \le f^*\fl_i \le {\bf 1}$. \item[(ii$'$)] We have $f^*\fm_i \ge \fE'$ (in $\bQ^{\Lam'}$) if and only if $i=t$.
\end{enumerate}
Hence by Theorem \ref{thm4-1} for $X'$, we have
\begin{align*}
\Ub{(f^\qc\fm_i)}M^q_{1,X'}\big/\Ub{(f^\qc\fm_{i+1})}M^q_{1,X'} & \cong \omega_{f^\qc\fm_i,f^\qc\fl_i,X'}^{q-1}\big/\cZ_{f^\qc\fm_i,f^\qc\fl_i,X'}^{q-1} \; \hbox{ for \; $0 \le i \le t-1$,} \\
\Ub{(f^\qc\fm_t)}M^q_{1,X'} & =0,
\end{align*}
and we are reduced to showing that
\begin{equation}\label{eq6-3}
R^1g_*\big(\omega_{f^{\qc}\fm,f^{\qc}\fl,X'}^{q-1}\big/\cZ_{f^{\qc}\fm,f^{\qc}\fl,X'}^{q-1}\big) = 0
\end{equation}
for $\fm,\fl \in \bN^{\Lam}$ with $\fm \ge 1$ and $\fl=\fd_\mu$ for some $\mu \in \Lam$.
We prove \eqref{eq6-3}. By \eqref{eq6-0}, there is a short exact sequence on $Y'_\et$
\[ 0 \lra \fcI^{(\fm+\fl)}\otimes_{\cO_{X'}}\tom_{Y'}^{q-1} \lra
 \fcI^{(\fm)}\otimes_{\cO_{X'}}\tom_{Y'}^{q-1} \lra \omega_{f^\qc\fm,f^\qc\fl,X'}^{q-1} \lra 0. \]
Since $\cI^{(\fm+\fl)}$ and $\cI^{(\fm)}$ are locally free $\cO_X$-modules, we obtain
\[ \omega_{\fm,\fl}^{q-1} \isom Rg_*\omega_{f^{\qc}\fm,f^{\qc}\fl,X'}^{q-1} \]
by applying $Rf_*$ to the above short exact sequence and projection formula together with the claim (a) in the proof of Lemma \ref{lem6-1}\,(1).
In particular we obtain
$R^1g_*\omega_{f^{\qc}\fm,f^{\qc}\fl,X'}^{q-1} = 0$ and \eqref{eq6-3}, because $R^1g_*$ is right exact for $p$-torsion sheaves for the reason of the dimension of the fibers of $g$ (\cite{SGA4} X.5.2, XII.5.2).
Thus we obtain \eqref{eq6-1}.
\par
By \eqref{eq6-1}, it remains to check
\begin{equation}\label{eq6-2}
R^1g_*\big(\Ua M_{1,X'}^q\big/\Ub{(f^{\qc}{\bf 1})}M_{1,X'}^q\big)=0
\end{equation}
If $a<c$, then we have $f^*{\bf 1} = {\bf 1}$ in $\bN^{\Lam'}$ and the assertion is obvious.
As for the case $c \le a$, we have
\[\Ua M_{1,X'}^q\big/\Ub{(f^{\qc}{\bf 1})}M_{1,X'}^q \cong \omega_{{\bf 1},\fd_\infty,X'}^{q-1}\big/\cZ_{{\bf 1},\fd_\infty,X'}^{q-1} \]
by \eqref{eq6++0} and Theorem \ref{thm4-1}.
We have
\[ Rg_* \omega_{{\bf 1},\fd_\infty,X'}^{q-1} \cong Rg_* \big(\cI_{X'}^{({\bf 1})}\big/\cI_{X'}^{(f^{\qc}{\bf 1})} \otimes_{\cO_{Y'}} g^*\tom_Y^{q-1} \big) \cong Rg_* \big(\cI_{X'}^{({\bf 1})}\big/\cI_{X'}^{(f^{\qc}{\bf 1})}\big) \otimes_{\cO_Y} \tom_Y^{q-1} \os{}=0 \]
by \eqref{eq6+1}.
Thus we obtain $R^1g_*\omega_{{\bf 1},\fd_\infty,X'}^{q-1} = 0$ and \eqref{eq6-2},
noting that $R^1g_*$ is right exact for $p$-torsion sheaves.
This completes the proof of Lemmas \ref{lem6-1} and \ref{lem5-3}, and Theorem \ref{thm5-1}.
\end{pf*}

\begin{rem}\label{cor6-2}
{\rm
By Lemma \ref{lem5-3} and the proof of Theorem \ref{thm5-1} (see the last diagram of \S\ref{sect5}),
we have $\Ua M_1^q \cong g_*(\Ua M_{1,X'}^q)$ as sheaves on $Y_\et$ for any $q \ge 0$ (under the notation of this section).
Consequently, we have
\[ \Ua M_1^q \isom  Rg_*(\Ua M_{1,X'}^q)
\quad \hbox{ in  \;\; $D(Y_\et)$ } \]
by Lemma {\rm \ref{lem6-1}\,(2)} and the fact that $R^s g_*(\Ua M_{1,X'}^q)=0$ for $s \ge2$, cf.\ \cite{SGA4} X.5.2, XII.5.2.}
\end{rem}


\subsection*{Acknowledgments}
The authors express their gratitude to Professor Kazuya Kato.
This paper is based on his notes on the computations on the complexes $\omega_{\fm,\fl}^\bullet$ (\S\ref{sect3}) and the filtration $\Ub {(\fm)}M^q_1$ (\S\ref{sect4}).
Thanks are also due to the referee, who read the manuscript carefully and gave the authors numerous constructive comments.


\end{document}